\documentclass[subeqn, letterpaper, reqno, 11pt]{amsart}
\usepackage[osf]{mathpazo}
\usepackage[T1]{fontenc}
\usepackage{amssymb, amsmath, amsthm, graphicx, mathtools, dsfont}
\usepackage[numbers]{natbib}
\usepackage[all]{xy}
\usepackage{hyperref}

\newcommand{\R}{\mathbb{R}}
\newcommand{\N}{\mathbb{N}}
\newcommand{\G}{\textsc{g}}
\newcommand{\E}{\mathbb{E}}

\newcommand{\Par}{\mathcal{P}}
\newcommand{\Hol}{\mathcal{H}}

\newcommand{\X}{\frak{X}}
\newcommand{\inv}{^{\text{-}1}}

\newcommand{\diam}{\textup{diam}}
\newcommand{\ptGH}{\xrightarrow{{}_{pt-GH}}}

\newcommand{\CvxRad}{\textup{CvxRad}}
\newcommand{\HolRad}{\textup{HolRad}}

\numberwithin{equation}{section}

\newtheorem{theorem}{Theorem}[section]
\newtheorem{lem}[theorem]{Lemma}
\newtheorem*{thma}{Theorem A}
\newtheorem*{thmb}{Theorem B}
\newtheorem*{thmc}{Theorem C}
\newtheorem*{thmd}{Theorem D}
\newtheorem*{thme}{Theorem E}

\newtheorem*{ndf}{Definition}

\newtheorem{proposition}[theorem]{Proposition}
\newtheorem{cor}[theorem]{Corollary}
\newtheorem{co}[theorem]{Conjecture}

\newtheorem{pro}{Problem}
\newtheorem{que}{Question}

\theoremstyle{definition}
\newtheorem{df}[theorem]{Definition}
\newtheorem{exa}[theorem]{Example}

\theoremstyle{remark}
\newtheorem{rem}[theorem]{Remark}
\newtheorem*{cla}{Claim}

\newcommand{\bd}{\begin{df}}
\newcommand{\ed}{\end{df}}

\newcommand{\bq}{\begin{que}}
\newcommand{\eq}{\end{que}}

\newcommand{\bl}{\begin{lem}}
\newcommand{\el}{\end{lem}}

\newcommand{\br}{\begin{rem}}
\newcommand{\er}{\end{rem}}

\newcommand{\bt}{\begin{theorem}}
\newcommand{\et}{\end{theorem}}

\newcommand{\bc}{\begin{cor}}
\newcommand{\ec}{\end{cor}}

\newcommand{\bco}{\begin{co}}
\newcommand{\eco}{\end{co}}

\newcommand{\bcl}{\begin{cla}}
\newcommand{\ecl}{\end{cla}}

\newcommand{\bp}{\begin{proposition}}
\newcommand{\ep}{\end{proposition}}

\newcommand{\brm}{\begin{rem}}
\newcommand{\erm}{\end{rem}}

\newcommand{\be}{\begin{equation}}
\newcommand{\ee}{\end{equation}}
\newcommand{\bx}{\begin{exa}}
\newcommand{\ex}{\end{exa}}
\newcommand{\bpr}{\begin{pro}}
\newcommand{\epr}{\end{pro}}

\begin{document}

\title{Convergence of Vector bundles with metrics of Sasaki-type}
\author{Pedro Sol\'orzano}
\address{State University of New York at Stony Brook}
\email{pedro@math.sunysb.edu}
\date{\today}
\dedicatory{To the memory of Detlef Gromoll.}

\maketitle
\begin{abstract} If a sequence of Riemannian manifolds, $X_i$, converges in the pointed Gromov-Hausdorff sense to a limit space, $X_\infty$,  and if $E_i$ are vector bundles over $X_i$ endowed with metrics of Sasaki-type with a uniform upper bound on rank, then a subsequence of the $E_i$ converges in the pointed Gromov-Hausdorff sense to a metric space, $E_\infty$.  The projection maps $\pi_i$ converge to a limit submetry $\pi_\infty$ and the fibers converge to its fibers; the latter may no longer be vector spaces but are homeomorphic to $\R^k/G$, where $G$ is a closed subgroup of $O(k)$ ---called the {\em wane group}--- that depends on the basepoint and that is defined using the holonomy groups on the vector bundles.  The norms $\mu_i=\|\cdot\|_i$ converges to a map $\mu_{\infty}$  compatible with the re-scaling in $\R^k/G$ and the $\R$-action on $E_i$ converges to an $\R-$action on $E_{\infty}$ compatible with the limiting norm. 
 
In the special case when the sequence of vector bundles has a uniform lower  bound on holonomy radius (as in a sequence of collapsing flat tori to a circle), the limit fibers are vector spaces.   Under the opposite extreme, e.g.  when a single compact $n$-dimensional manifold is re-scaled to a point, the limit fiber is $\R^n/H$ where $H$ is the closure of the holonomy group of the compact manifold considered.

An appropriate notion of parallelism is given to the limiting spaces by considering curves whose length is unchanged under the projection. The class of such curves is invariant under the $\R$-action and each such curve preserves norms. The existence of parallel translation along rectifiable curves with arbitrary initial conditions is also exhibited. Uniqueness is not true in general, but a necessary condition is given in terms of the aforementioned  wane groups $G$. 
\end{abstract}

\section*{Introduction}
The Gromov-Hausdorff convergence of Riemannian manifolds was introduced by Gromov in the late 1970's.  Unlike smooth convergence, limit spaces under Gromov-Hausdorff convergence need not be smooth or even Lipschitz.  Adding conditions of uniform curvature bounds to the sequence of metrics one can control the regularity of the limit spaces to some extent.  Work of  \citet{MR1484888, MR1815410, MR1815411} showed that the regular set of limit spaces with Ricci curvature bounded below is dense and a $C^{1.\alpha}$ submanifold.   Adding the stronger condition of one- or two-sided bounds on sectional curvature, there have been significantly stronger results, many structural results have been obtained by \citet{MR0263092}, \citet*{MR1126118}, \citet{MR1097241}, \citet{MR1863020}, as well as upcoming work of  \citet{RONG}. Only by assuming both a sectional curvature and lower bound on volume or that the sequence is Einstein with a lower bound on injectivity radius does one obtain limits which are $C^{1.\alpha}$ manifolds, as seen by \citet{MR1074481}, and \citet{MR1158336}. However, it should be noted that a common feature is to have certain assumptions on the curvature, injectivity radius, etc. 

In this communication it is only assumed that there is a sequence of Riemannian manifolds which converges in the pointed Gromov-Hausdorff sen\-se to a limit space; the sequence of tangent bundles over those Riemannian manifolds ---or more generally, an arbitrary sequence of vector bundles over the converging sequence of Riemannian manifolds---  is then analyzed. Throughout this communication, none of the usual uniform bounds on curvature, diameter, volume or injectivity radius are assumed.   Only those properties which can be derived from the pointed Gromov-Hausdorff convergence of the base spaces are used.

It is worth noting that the results discussed here differ from the very interesting approach taken by \citet{MR2600284}. He introduces a Lipschitz seminorm of a very natural space of matrix-valued functions to control distances between vector bundles.  In essence, he regards Euclidean vector bundles as a certain type of map into the space of self-adjoint idempotent matrices. In the case of vector bundles over smooth manifolds, this can be easily be seen as maps into a suitable Grassmannian.  Under the assumption that two (compact) metric spaces be $\varepsilon$-close, he gives a correspondence between their vector bundles with control on their Lipschitz seminorm on any metric on their disjoint union that  makes said spaces $\varepsilon$-Hausdorff close.

Vector bundles with metric connections (i.e. Euclidean bundles with compatible connections) have natural metrics Riemannian metrics on their total spaces called metrics of Sasaki-type (see Definition \ref{sasmetdef}).  These metrics were first introduced on tangent bundles by \citet{MR0112152} and for more general vector bundles by \citet*{0703059}, where they introduce a two-parameter family of metrics that include a metric know to exist by the results of \citet{MR0309010}. These metrics coincide with the full classification of the natural metrics on tangent bundles given by \citet{MR974641}. It should be noted that the use of the word {\em natural} coincides with its usage in the classification of natural bundles given by \citet{MR509074} as part of her doctoral dissertation; and thus the results stated here could be stated as certain continuity properties of these {\em natural} bundle functors.

The first explicit rendering of the Cheeger-Gromoll metric, together with a systematic study of the Sasaki metric was given by \citet{MR946027}. Later, it has been developed by many, in particular by \citet{MR2106375}. Most of the attention has been for the case of the tangent bundle.  In this case, on the tangent bundle, $TM$, over a Riemannian manifold, $M$, with the standard connection, the Sasaki metric on $TM$ is uniquely defined so that $\pi:TM \to M$ be a Riemannian submersion where the horizontal lifts of curves are simply parallel translations along curves and, furthermore, that the individual tangent spaces be totally geodesic flats; i.e with the intrinsic distance, the fibers are isometric to Euclidean space.  However, with the restricted metric, distances between points in a fiber may be achieved by paths that leave the fibers; in some cases even by horizontal paths, thus relating the problem to the semi-Riemannian context.  The fibers with the restricted metrics are holonomic spaces \cite{1004.1609}, whose metrics depends on the holonomy group and the shortest lengths of curves representing each holonomy element (See Definition \ref{distsas}).   

In Example  \ref{totcollapse}, if a single compact $n$-dimensional Riemannian manifold 
is re-scaled so that it converges in the Gromov-Hausdorff sense to a single point, then the Gromov-Hausdorff limit of the tangent bundles endowed with their Sasaki metrics converge to $\R^n/H$ where $H$ is the closure of the holonomy group of $M$ inside the orthogonal group $O(n)$.  Intuitively, horizontal paths became so short under rescaling that in the limit, vectors related by a horizontal curve are no longer distinct points. In contrast, if one has a sequence of standard 2 dimensional flat tori collapsing to a circle (Example \ref{prodcollapse}), then the limit of the tangent bundles is  $S^1 \times \E^2$, where the fibers are Euclidean since the holonomy group is trivial in this setting.  

Observe also that there are sequences of Riemannian manifolds which converge in the Gromov-Hausdorff sense whose tangent bundles do not converge (Example   \ref{prodcollapse}).  
Nevertheless a precompactness theorem for tangent bundles and other vector bundles can be obtained:

\begin{thma}[Theorem \ref{precptBWC}] Given a precompact collection of (pointed) Riemannian manifolds $\mathcal{M}$ and a positive integer $k$, the collection $\textup{BWC}_k(\mathcal{M})$ of vector bundles with metric connections of rank $\leq k$ endowed with metrics of Sasaki-type is also precompact. The distinguished point for each such bundle is the zero section over the distinguished point of their base.
\end{thma}

The assumption that the rank be bounded is easily satisfied for natural bundles over a convergent sequence of Riemannian manifolds (such as tangent bundles, cotangent bundles, or combinations thereof). 

In \cite{1004.1609}, the author has investigated the restricted metric of the fibers of a bundle with metric and connection endowed with its associated metric of Sasaki-type. In Section \ref{holspsec} the notion of holonomic space metric is review and in Section \ref{GHhol} the convergence of these metrics is analyzed to try to understand the fiberwise behavior of a convergent sequence of metrics of Sasaki-type. This approach proved to be quite useful in view of the following results.

\begin{thmb}[Propositions \ref{sascpt}, \ref{fibsGH}, \ref{equidfib} and Corollary \ref{absres}] For any sequence of Riemannian manifolds $\{(X_i,p_i)\}$ converging to $(X_{\infty},x_{\infty})$ consider a convergent family of bundles with metric connection $(E_i,h_i,\nabla_i)$ over it converging to $(E_{\infty},y_{\infty})$. Then there exist continuous maps $\pi_{\infty}:E_{\infty}\rightarrow X_{\infty}$, $\varsigma_{\infty}:X_{\infty}\rightarrow E_{\infty}$, $\mu_{\infty}:E_{\infty}\rightarrow\R$, and a subsequence, without loss of generality also indexed by $i$, such that:
\begin{enumerate}
\item the projection maps $\pi_i:E_i\rightarrow X_i$ converge to $\pi_{\infty}:E_{\infty}\rightarrow X_{\infty}$, which is also a submetry with equidistant fibers;
\item the zero section maps $\varsigma_i:X_i\rightarrow E_i$ converge to $\varsigma_{\infty}:X_{\infty}\rightarrow E_{\infty}$, which is also a isometric embedding; 
\item $\pi_{\infty}\circ \varsigma_{\infty}=id_{X_{\infty}}$; 
\item the maps $\mu_i:E_i\rightarrow\R$, given by 
$$\mu_i(u)=d_{E_i}(u,\varsigma_i\circ\pi_i(u))=\sqrt{h_i(u,u)},$$ 
converge to $\mu_{\infty}:E_{\infty}\rightarrow\R_{\geq0}$ also given by $$\mu_{\infty}(y)=d_{E_{\infty}}(y,\varsigma_{\infty}\circ\pi_{\infty}(y));$$
\item The scalar multiplications on $E_i$ converge to an $\R-$action on $E_{\infty}$ such that
$$\mu_{\infty}(\lambda u)=|\lambda|\mu_{\infty}(u)$$
\end{enumerate}
For any $\varepsilon>0$ and for any sequence $\{q_i\}$, $q_i\in X_i$, converging to $q\in X_{\infty}$,
\begin{enumerate}
\setcounter{enumi}{5}
\item $\pi_i\inv(B_{\varepsilon}(q_i))\ptGH\pi_{\infty}\inv(B_{\varepsilon}(q))$;
\item $\pi_i\inv(q_i)\ptGH\pi_{\infty}\inv(q)$.
\end{enumerate}
\end{thmb}

As mentioned before, Example \ref{totcollapse} already suggests that the holonomy group must play a significant r\^ole. In the non-symmetric simply-connected setting, Berger gave a complete classification of the possible holonomy representations for the Levi-Civita connection \cite{MR0079806}. In particular, there are finite\-ly many, once an upper bound in dimension is assumed. With this in mind, the following result yields more information about fibers of the limiting map.

\begin{thmc}[Theorem \ref{finholtype}] Let $\pi_i:E_i\rightarrow X_i$ be a convergent sequence of vector bundles with bundle metric and compatible connections $\{(E_i,h_i, \nabla_i)\}$, with limit $\pi:E\rightarrow X$. Suppose further that there are only finitely many holonomy representation types. Then there exists a positive integer $k$ such that for any point $p\in X$ there exists a compact Lie group $G\leq O(k)$ , that depends on the point, such that the fiber $\pi\inv(p)$ is homeomorphic to $\R^k/G$, i.e. the orbit space under the standard action of $G$ on $\R^k$.
\end{thmc}

The group $G$ here is described explicitly in Theorem \ref{dinfty}  and will be called the {\em wane group} at $x\in X$ because of another precise description of
the fibers as $V/G_0$ given in Theorem \ref{limsup} where $G_0$ is defined in terms of the metrics $d_i$  (on the fibers of the converging sequence of vector bundles) by essentially looking at sequences of holonomy elements with waning norm.

It is important to remark that the wane group G truly depends on the base point; thus the limit $\pi: E\to X$ need not be a fiber bundle.  One can imagine this occurs for a example when a sequence of Riemannian manifolds converge smoothly everywhere except at a point and at the point they
develop a conical singularity.  In that setting one expects the fibers away from the singularity
to be vector spaces while the singular point's fiber is not Euclidean.  

In \cite{1004.1609}, the author introduced the notion of {\em holonomy radius at a point}. Because for metrics of Sasaki-type the fibers of the vector bundle in cosideration are totally geodesic and flat, it makes sense to consider the following definition.

\begin{ndf} Consider a vector bundle $E$ with metric and connection over a Riemannian manifold.  The {\em holonomy radius} of a point $p$ in the base is the largest $R>0$ such that the restricted metric on $B_R(0_p)\cap E_p\subseteq E$ is Euclidean.
\end{ndf}
For a more technical definition see Definition \ref{holraddef}.  Now, if one is willing to assume some further restriction on a convergent sequence of manifolds, a uniform lower bound on their holonomy radii yields the following.

\begin{thmd}[Theorem \ref{hlradbnd}] Let $\pi_i:E_i\rightarrow X_i$ be a convergent sequence of vector bundles with bundle metric and compatible connections $\{(E_i,h_i, \nabla_i)\}$, with limit $\pi:E\rightarrow X$. Suppose further that there are only finitely many holonomy representation types  and that there exist a uniform positive lower bound for the holonomy radii of $\pi_i:E_i\rightarrow X_i$. Then the fibers of  $\pi_{\infty}$ are vector spaces. 
\end{thmd}

In the context of Alexandrov spaces, the notion of tangent cone is closely related to the usual notion of tangent vector in the sense of realizing it as a cone over the space of direction, as noted by \citet*{MR1185284}. It has yet to be studied how this relates to the fibers of the limits of tangent bundles when the bases have lower curvature bounds. The fibers are seen to be quotients of vector spaces by groups that act on spheres; as such, they are certainly topologically cones.

Furthermore, there is a natural way to define a notion of parallelism on these limit spaces by considering horizontal curves.

\begin{ndf} Given a submetry $\pi: Y\rightarrow X$, a curve in $Y$ is horizontal if and only if its length is equal to the length of its projection in $X$
\end{ndf}

The collection of horizontal curves over a given curve $\alpha$ in $X$ gives a relation between the fibers the endpoints of $\alpha$. For loops at a point, it follows that the set of parallel translates form a ${}^*$-semigroup, which will be called the {\em Holonomy monoid} of $\pi$, because it generalizes the holonomy group, yet it is not necessarily a group.

In the particular of the limit spaces $\pi: E\rightarrow X$ satisfy very nice properties summarized in the next result.  

\begin{thme}[Corollary \ref{absnorm}] Let $\pi_i:E_i\rightarrow X_i$ be a convergent sequence of vector bundles with bundle metric and compatible connections $\{(E_i,h_i, \nabla_i)\}$, with limit $\pi:E\rightarrow X$. Given any curve $\alpha:I\rightarrow X$ and a point $u\in\pi\inv(\alpha(0))$ there exists a parallel translate $\gamma$ of $\alpha$ with initial point $u$, furthermore the norm is constant along $\gamma$ and any re-scaling of $\gamma$ is also a parallel translate of $\alpha$.
\end{thme}

The non uniqueness of parallel translates is exactly encoded by the lack of invertibility of holonomy elements (see Theorem \ref{holmonpar}). Also, a necessary condition for having uniqueness for this weak notion of parallelism is given by whether the wane groups are conjugates of each other or not (see Theorem \ref{parwane}).

This notion of parallelism should be compared to that of \citet{MR1601854}, again in the case of lower curvature bounds.

\vspace{11pt}

{\bf\noindent Acknowledgments.} This communication is dedicated to the memory of Detlef Gromoll, to whom the author is greatly indebted; for all his teachings, his support, and for his contagious love for mathematics and life. The author is also very thankful for the complete support and the (perhaps undeserved) faith he has received from his doctoral advisers, H. Blaine Lawson Jr. and Christina Sormani. 

The author also wishes to thank Kristopher Tapp and Xiaochun Rong for their insight and for their useful remarks.


\setcounter{tocdepth}{3}
\tableofcontents


\section{Background} Several distinct approaches to the understanding of Riemannian manifolds will interplay in the sequel. In particular, the geometric properties induced on the total space of their tangent bundle. It should be noted that here {\em geometric} has two different meanings. One is the usual differential-geometric approach, where a smooth Riemannian metric (of Sasaki-type) is induced to the total space of vector bundles over Riemannian manifolds and several ---quite rigid--- properties have been obtained (see \cite{MR0286028, MR946027}). Another path is the metric geometric approach; i.e. the study of the metric spaces obtained from the infinitesimal metrics. The fibers have two distinct metrics on them inherited from that of the total space. One is known to be Euclidean, the other was observed by the author to be of a certain very explicit type (called {\em holonomic metric} in \cite{1004.1609}). 

Finally, the degeneration of these structures investigated in the communication is mainly through Gromov's extended notion of the Hausdorff distance; which has proved to be quite a fecund area of research in recent years.  

\subsection{Gromov-Hausdorff convergence}
Listed are a few results needed in the sequel. For a detailed treatment of this concepts see \cite{MR1835418,MR2307192, MR2243772}. 

\bd[\citet{MR2307192}] Given two complete metric spaces $(X,d_X)$ and $(Y,d_Y)$, their {\em Gromov-Hausdorff} distance is defined as the following infimum.
\be\label{GHdist}
d_{GH}(X,Y)=\inf\left\{\varepsilon>0\left|\begin{array}{rl}(1)&\exists d:(X\sqcup Y)\times(X\sqcup Y)\rightarrow\R,\text{ metric}\\ (2)& d|_{X\times X}=d_X,d|_{Y\times Y}=d_Y \\ (3)&\forall x\in X(\exists y\in Y, d(x,y)<\varepsilon) \\ (4)&\forall y\in Y(\exists x\in X, d(x,y)<\varepsilon)\end{array}\right. \right\}
\ee

That is the infimum of possible $\varepsilon>0$ for which there exists a metric on the disjoint union $X\sqcup Y$ that extends the metrics on $X$ and $Y$, in such a way that any point of $X$ is $\varepsilon$-close to some point of $Y$ and vice versa.  
\ed
\br This is a generalization of the {\em Hausdorff distance} between subspaces of a fixed metric space $(Z,d)$. In this case, the distance $d^Z_H(X,Y)$, between subspaces $X,Y\subseteq Z$, is defined as follows.
\be\label{Hdist}
d_{H}(X,Y)=\inf\left\{\varepsilon>0\left|\begin{array}{rl}(3)&\forall x\in X(\exists y\in Y, d(x,y)<\varepsilon) \\ (4)&\forall y\in Y(\exists x\in X, d(x,y)<\varepsilon)\end{array}\right. \right\}
\ee

\er
\br For compact metric spaces the assignment \eqref{GHdist} is always finite, since 
$$d_{GH}(X,Y)\leq \frac{1}{2}\max\{diam(X),diam(Y)\};$$
it may however be infinite if compactness is not assumed. This assignment is positive, symmetric and satisfies the triangle inequality (provided it makes sense). Two (compact) spaces are zero distance apart if and only if they are isometric. 
\er
\bd[\citet{MR2307192}] Let $X$ and $Y$ be metric spaces. For $\varepsilon>0$, an {\em $\varepsilon$-isometry} from $X$ to $Y$ is a (possibly non-continuous) function $f:X\rightarrow Y$ such that:
\begin{enumerate}
\item for all $x_1,x_2\in X$,
\be
\left|d_X(x_1,x_2)-d_Y(f(x_1),f(x_2))\right|<\varepsilon;\text{ and}
\ee 
\item for all $y\in Y$ there exists $x\in X$ such that
\be 
d_Y(f(x),y)<\varepsilon.
\ee
\end{enumerate}
\ed

\bp[\citet{MR2307192}]\label{epsiso} Let $X$ and $Y$ be metric spaces and $\varepsilon>0$. Then,
\begin{enumerate}
\item if $d_{GH}(X,Y)<\varepsilon$ then there exists a $2\varepsilon$-isometry between them. 
\item if there exists an $\varepsilon$-isometry form $X$ to $Y$, then $d_{GH}(X,Y)<2\varepsilon$.
\end{enumerate}
\ep

Except for some set theoretical considerations, the collection of isometry classes of metric spaces, together with the Gromov-Hausdorff distance, $(\mathcal{M},d_{GH})$ behaves like an extended metric space (i.e. allowing infinite values). When restricted to compact metric spaces, it is a metric space and, as such, yields a notion of convergence for sequences.

\br\label{GHisH} It was proved by \citet{MR2307192} that if a sequence $\{X_i\}$ of compact metric spaces converges in the Gromov-Hausdorff sense to a compact metric space $X$, then there exists a metric on $X\sqcup\bigsqcup_iX_i$ for which the sequence $\{X_i\}$ converges in the Hausdorff sense. Because of this, it makes sense to say that a sequence of points $x_i\in X_i$ converge to a point $x\in X$. 

\er

\br\label{ptsqnHaus} In this setting, a sequence of subspaces $X_i\subseteq Z$ converges to a subspace $X\subseteq Z$ if and only if: 
\begin{enumerate}
\item for any convergent sequence $x_i\rightarrow x$, such that $x_i\in X_i$ for all $i$, it follows that $x\in X$; and
\item for any $x\in X$ there exists a convergent sequence $x_i\rightarrow x$, with $x_i\in X_i$.
\end{enumerate}

\er

\bd[\citet{MR2307192}]\label{ptsqns} Let $\{X_i\}$, $\{Y_i\}$ be convergent sequences of pointed metric spaces and let $X$ and $Y$ be their corresponding limits. One says that a sequence of continuous functions $\{f_i\}:\{X_i\}\rightarrow\{Y_i\}$ {\em converges} to a function $f:X\rightarrow Y$ if there exists a metric on $X\sqcup\bigsqcup_iX_i$ for which the subspaces $X_i$ converge in the Hausdorff sense to $X$ and such that for any sequence $\{x_i\in X_i\}$ that converges to a point $x\in X$, the following holds.
\be 
f(x)=\lim_{i\rightarrow\infty}f_i(x_i)
\ee
\ed
\br The limit function $f$ is unique if it exists; i.e. it is independent of the choice of metric on $X\sqcup\bigsqcup_iX_i$.
\er

The following is Gromov's way to produce a notion of convergence for the non-compact case. For technical reasons, the assumption that the spaces be proper (i.e. that the distance function from a point is proper, thus yielding that closed metric balls are compact) is required \cite{MR2307192}. 

\bd A sequence $\{(X_i,x_i)\}$ of pointed proper metric spaces is said to converge to $(X,x)$ in the pointed Gromov-Hausdorff sense if the following holds: For all $R>0$ and for all $\varepsilon>0$ there exists $N$ such that for all $i>N$ there exists an $\varepsilon$-isometry 
$$f_i:B_R(x_i)\rightarrow B_R(x),$$
 with $f_i(x_i)=x$, where the balls are endowed with restricted (not induced) metrics. 
\ed
\br In the previous definition, it is enough to verify the convergence on a sequence of balls around $\{x_i\}$ such that their radii $\{R_j\}$ go to infinity. Furthermore, the limit is necessarily also proper as noted by \citet{MR623534}. 
\er

\br
Given a pointed space $(X,x)$, \citet{MR2307192} studies the relation between precompactness and the function that assigns to each choice of $R>0$ and $\varepsilon>0$ the maximum number $N=N(\varepsilon, R,X)$ of disjoint balls of radius $\varepsilon$ that fit within the ball of radius $R$ centered at an $x \in X$. Furthermore he proves the following result.
\er

\bt[Gromov's Compactness Theorem \cite{MR2307192}, Prop.5.2]\label{Gromovs} A family of pointed path metric spaces $(X_i,x_i)$ is pre-compact with respect to the pointed Gromov-Hausdorff convergence if and only if each function $N(\varepsilon, R, \cdot)$ is bounded on $\{X_i\}$. In this case, the family is relatively compact, i.e., each sequence in the $X_i$ admits  a subsequence that converges in the pointed Gromov-Hausdorff sense to a complete, proper path metric space. 
\et

\br \label{GromovN}Providing a bound for $N$ is equivalent to providing a bound to the minimum number of balls of radius $2\varepsilon$ required to cover the ball of radius $R$ (see \cite{MR2243772}). This will be used instead in the sequel and will also be denoted by $N$.
\er

\br In the non-compact setting, in order to consider the convergence of sequences of points $\{p_i\in X_i\}$, as in Remark \ref{GHisH}, the only technicality is the following: In order for a sequence $\{p_i\}$ to be convergent, it has to be bounded. Therefore, there must exist a large enough $R>0$ such that for all $i$, $p_i\in B_R(x_i)$, where the $x_i\in X_i$ are the distinguished points. Because of this, a sequence  $\{p_i\in X_i\}$ is convergent if there exists $R>0$ for which the sequence  $\{p_i\in \overline{B_R(x_i)}\subseteq X_i\}$ is convergent as in Remark \ref{GHisH}.
\er

To analyze the behavior of sequences of functions defined on convergent sequences of spaces, \citet{MR623534}, as well as \citet{MR1094462}, has given a generalization to the classical Arzel\`a-Ascoli Theorem. Their setting is that of compact spaces. To analyze the non-compact setting, a further generalization is required.  In the proof of the compact case, families of countable dense sets $A_i\subseteq X_i$, $A\subset X$ are considered (since continuous functions are determined by their values on dense sets, and as a basis for the standard diagonalization argument). Also, the codomains satisfy that for every sequence there exist a convergent subsequence. To retain these properties, the assumption of separability for the domains and the requirement of totally bounded metric balls for the codomains are added; these are both controlled by the assumption that all the spaces being considered be proper. Also, by virtue of the Hopf-Rinow Theorem in Riemannian geometry, when the metric spaces considered are Riemannian manifolds these conditions follow from completeness.  

\bt[Arzel\`a-Ascoli Theorem, cf. \cite{MR623534, MR1094462}]\label{ArAs}
Consider two convergent sequences of complete proper pointed metric spaces, $\{(X_i,x_i)\}$, $\{(Y_i,y_i)\}$. Let $(X,x)$ and $(Y,y)$ be their corresponding limits. Suppose further that there is an equicontinuous sequence of continuous maps $\{f_i\}$,
\be 
f_i:X_i\rightarrow Y_i,
\ee
such that $f_i(x_i)=y_i$, for all $i$. Then there exist a continuous function $f:X\rightarrow Y$, with $f(x)=y$, and a subsequence of $\{f_i\}$ that converges to $f$. 
\et
\begin{proof} In the compact case, this is the content of the generalization of the Arzel\`a-Ascoli Theorem given by \citet{MR1094462}. For the non-compact case, \citet{MR623534} already proved this for isometries under the assumption of properness. Again, a diagonalization argument is required. Namely, because the sequence $\{f_i\}$ is equicontinuous, for every $\varepsilon>0$ consider the largest $\delta=\delta(\varepsilon)$ that satisfies the definition of equicontinuity. It follows that $\delta$ is an increasing function of $\varepsilon$ that goes to infinity as $\varepsilon$ does; it may happen that $\delta(\varepsilon)=\infty$ for a finite $\varepsilon$. This implies that for any $R>0$ there exists $\tilde R>0$,
\be
f_i(B_R(x_i))\subseteq B_{\tilde R}(y_i),
\ee
by essentially considering the inverse of $\delta$ as a function of $\varepsilon$ (if $\delta$ is infinite, then the existence of $\tilde R$ is clearly also satisfied.). This means that one can now repeat the proof of the compact case for the restrictions $\{f_i|_{B_R(x_i)}\}$ (since $B_R(x_i)$ is separable and $B_{\tilde R}(y_i)$ compact). Therefore, consider a sequence of radii $R_j\rightarrow\infty$ and apply the standard diagonalization argument to the successive restrictions of the convergent subsequences of $f_i$'s  (and subsequences thereof) to $B_{R_j}(x_i)\rightarrow B_{\tilde R_i}(y_i)$. By uniqueness of the limit, one obtains further and further extensions to a single continuous function $f: X\rightarrow Y$ as promised.
\end{proof}

\subsection{Metrics of Sasaki-type} \label{Sasmetsec}
Given a Riemannian manifold $(M,g)$, from all the metrics that render the total space of the tangent bundle a Riemannian manifold itself, and the standard projection a Riemannian submersion, perhaps the most natural pick is the metric introduced by \citet{MR0112152}. In a more general setting, a Sasaki-type metric can be introduced to any bundle over a Riemannian manifold, provided that the bundle be further equipped with a metric connection. 
For a more detailed study of the aforementioned topics see \cite{1004.1609}, as well as the book by \citet{MR1393940}.

\bd[Connections and metrics on vector bundles]   Given a vector bundle $\pi:E\rightarrow M$ over a smooth manifold $M$, a bundle metric is a choice of inner products on each fiber, $E_p$, that depends smoothly on the base space. Namely, it is a smooth bundle map $h:E\oplus E\rightarrow \underline{\R}$, where $\underline{\R}$ is the trivial bundle $M\times\R$, and such that for any $p\in M$ the map $h_p$ is an inner product on $E_p$. Yet another interpretation comes from regarding $h$ as a section of $Sym^2(E^*)$.

A connection on a bundle $\pi:E\rightarrow M$ is a map $\nabla$, 
\be 
\nabla:\Gamma(E)\rightarrow\Gamma(T^*M\otimes E),
\ee
that satisfies the Leibniz rule $\nabla(f\sigma)=df\otimes\sigma+f\nabla(\sigma)$ for any section $\sigma$ and any smooth function $f$ on $M$. Given a bundle metric, the connection is said to be {\em metric} if it further satisfies that $\nabla(h)=0$, where $\nabla$ the induced connection on $\text{Sym}^2E^*$. More explicitly, a connection is metric if and only if for any sections $\sigma, \tau\in\Gamma(E)$ and any vector field $X\in\X(M)$, 
\be
X(h(\sigma,\tau))= h(\nabla_X\sigma, \tau)+h(\sigma,\nabla_X\tau).
\ee
\ed

\bd[Parallel translation]
Given a vector bundle $\pi:E\rightarrow M$ with connection $\nabla$, a section $\sigma\in\Gamma(E)$ is {\em parallel} if $\nabla\sigma\equiv0$. A section $\sigma$ along a curve $\alpha:[0,1]\rightarrow M$ is {\em parallel} if $\nabla_{\dot\alpha}\sigma\equiv0$.  Given any curve $\alpha$ and a vector $u\in E$ with $\pi(u)=\alpha(0)$ there exists a unique parallel section $t\mapsto P^{\alpha}_t(u)$ along $\alpha$ with $P^{\alpha}_0(u)=u$. 

It follows that the transformation $u\mapsto P^{\alpha}_t(u)$ is linear with respect to $u$ for any $t$. Furthermore, if the connection is metric then $P^{\alpha}_t(u)$ is an isometry with respect to the bundle metric. $P^{\alpha}_t(u)$ is frequently called {\em parallel translation} of $u$ along $\alpha$ at time $t$.
\ed

\bd[Holonomy Groups] Given a vector bundle $\pi:E\rightarrow M$ with connection $\nabla$, and given any $p\in M$, the {\em holonomy group} of $\nabla$ is the collection, denoted by $Hol_p(\nabla)$, of $P^{\alpha}_1:E_p\rightarrow E_p$ where $\alpha:[0,1]\rightarrow M$ is a loop at $p$; i.e. $\alpha(0)=\alpha(1)=p$. If $M$ is connected it follows that, for all point $p,q\in M$, $Hol_p(\nabla)$ is isomorphic to $Hol_q(\nabla)$, but this isomorphism is not canonical. Furthermore, if $\nabla$ is metric with respect to $h$, then for all $p\in P$, $Hol_p(\nabla)$ is a subgroup of the orthogonal group $O(E_p)$ with respect to $h_p$.
\ed

In order to define metrics of Sasaki-type, the following construction is required. 
\bd[Vertical lifts] 
Given a vector bundle $\pi:E\rightarrow M$, consider a vector $u\in E_p$, the {\em vertical lift} of $u$ is the map $u^v:E_p\rightarrow T(E_p)\subseteq TE$ given by
\be v\mapsto \dot\gamma(0),
\ee
where $\gamma$ is the curve in $E$ given by $\gamma(t)=v+tu$.  It follows that $\pi_*(u^v)\equiv0$.

In fact, if one starts with a section $\sigma\in \Gamma(E)$, in this fashion one produces a vector field $\sigma^v\in\X(E)$, the vertical lift of $\sigma$, that satisfies that 
\be 
\pi_*(\sigma^v)\equiv0.
\ee 
\ed

\bd[Horizontal lifts]  Given a vector bundle $\pi:E\rightarrow M$ with connection $\nabla$, consider a tangent vector $x\in T_pM$, the {\em horizontal lift} of $x$ is the map $x^h:E_p\rightarrow TE$ given by
\be
v\mapsto \dot\sigma(0),
\ee
where $\sigma(t)=P^{\alpha}_t(v)$ and $\alpha$ is any curve on $M$ such that $\alpha(0)=p$ and $\dot\alpha(0)=x$; i.e. $x^h(v)$ is the derivative of the parallel translation of $v$ in the direction of $x$. It follows that $\pi_*(x^h)\equiv x$.

In fact, if one starts with a vector field  $X\in \X(M)$, in this fashion one produces a vector field $X^h\in\X(E)$, the horizontal lift of $\sigma$, that satisfies that 
\be
\pi_*(X^h)\equiv X.
\ee 
\ed

\br Given a vector bundle $\pi:E\rightarrow M$ with metric connection $\nabla$ and bundle metric $h$, consider any vector $\xi\in T_uE$, with $p=\pi(u)$, then $\xi$ can be expressed as
\be
\xi=\sigma^v(u)+x^h(u)
\ee
for some uniquely determined $x\in T_pM$ and $\sigma\in E_p$. 
\er
\bd[Metrics of Sasaki-type]\label{sasmetdef} Given a vector bundle $\pi:E\rightarrow M$ over a Riemannian metric $(M,g)$, a bundle metric $h$ on $E$, and a metric connection $\nabla$, the {\em Sasaki-type metric} associated to these data, denoted by $\G=\G(g,h,\nabla)$, is given in the following way.  For any vector fields $X,Y\in\X(M)$ and for any sections $\sigma,\tau\in\Gamma(E)$:
\begin{align}
\G_u(\sigma^v,\tau^v)&=h_{\pi(u)}(\sigma,\tau),\\
\G_u(\sigma^v,Y^h)&=0,\\
\G_u(X^h,Y^h)&=g_{\pi(u)}(X,Y).
\end{align}
\br In the case when $E=TM$, the Levi-Civita connection $\nabla$ with respect to a Riemannian metric $g$ and $h=g$ already provide all the data required for this definition. This is then the Sasaki metric of the tangent bundle (cf. \citet{MR0112152}). 
\er
\ed

\br\label{sasprops} The differential geometric properties of these metrics have been extensively studied to induce metrics on unit tangent bundles, to consider harmonic vector fields, or even to study their naturality in the category of Riemannian maps, etc. (see \cite{MR2330461, MR1798731,MR974641,MR946027,1004.1609}). Some of its nice features are:
\begin{enumerate}
\item The Sasaki-type metric $\G$ is complete if and only if $g$ is also complete.
\item The projection $\pi:E\rightarrow M$ is a Riemannian submersion. 
\item The fibers $E_p$ are totally geodesic and flat.
\item The zero section $\varsigma:M\rightarrow E$, $\varsigma(p)=0_p\in E_p$ is an isometric embedding (i.e. it is distance preserving).
\item The rays $t\mapsto tu$ are geodesic rays for $t\geq0$.
\end{enumerate}
\er

\bp[\citet{1004.1609}]\label{distsas} The length distance on $(E,\G)$ is expressed as follows. Let $u,v\in E$
\be\label{distsaseq}
d_{E}(u,v)=\inf\left\{\sqrt{\ell(\alpha)^2+\|P^{\alpha}_1u-v \|^2} \left|\begin{matrix} \alpha:[0,1]\rightarrow M,\\ \alpha(0)=\pi(u),\\ \alpha(1)=\pi(v)\end{matrix}\right.\right\}.
\ee
Furthermore, if $\pi (u)=\pi (v)$ then
\be\label{distfor}
d_{E}(u,v)=\inf\{\sqrt{L(a)^2+\|au-v \|^2} : a\in Hol_p\},
\ee
with $L$ being the infimum of lengths of loops yielding a given holonomy element.
\ep

\subsection{Holonomic Spaces}\label{holspsec}

In \cite{1004.1609}, the author introduced the notion of holonomic space as a way to understand the geometry of the restricted metric of fibers of vector bundle with metric connection over a Riemannian manifold, when said bundles are endowed with their corresponding metrics of Sasaki-type (see Subsection \ref{Sasmetsec}).  Listed here are the basic properties of such spaces that will be needed in the sequel.

\bd[see \cite{NORMED}] Given an abstract group $H$ with identity element $e$, a {\em group-norm} on $H$ is a function $L:H\rightarrow\R$ that satisfies the following properties. For all $a,b\in H$,
\begin{enumerate}
\item $L(a)\geq0$ and $L(a)=0$ if and only if $a=e$;
\item $L(a\inv)=L(a)$; and
\item $L(ab)\leq L(a)+L(b)$.
\end{enumerate}
\ed

\bd[\citet{1004.1609}] A {\em holonomic space} is a triplet $(V,H,L)$: $V$ is a normed vector space, $H$ is a group of norm preserving linear maps, and $L$ is a group-norm on $H$; such that together satisfy that for any $u \in V$ there exists a positive number $R>0$ such that if $v,w\in V$ and for any $a\in H$, if
$$\|u-v\|,\|u-w\|<R$$
then
\be\label{propP}
\|v-w\|^2-\|av-w\|^2\leq L^2(a).
\ee
\ed

\bt[\citet{1004.1609}] \label{hlmetthm} Let $(V,H,L)$ be a holonomic space.  
\be\label{hlmet} d_L(u,v)=\inf_{a\in H}\left\{\sqrt{L^2(a)+\|u-av\|^2}\right\},
\ee
is a metric on $V$. Furthermore, the identity map is a distance non-increasing map between $(V,\|\cdot\|)$ and $(V,d_L)$.
\et

\bd[\cite{1004.1609}]  The metric given by \eqref{hlmet} is called {\em holonomic space metric}.
\ed

\bd[\cite{1004.1609}] \label{hlsprad} Given a holonomic space $(V,H,L)$, 
\begin{enumerate}
\item the {\em Holonomy radius} at $u\in V$, denoted by $\HolRad(u)$, is the supremum of the radii $R>0$ such that the metric $d_L$ is Euclidean when restricted to the ball of radius $R$; equivalently, the supremum of $R>0$ for which  \eqref{propP} holds.
\item the {\em Convexity radius} of $(V,H,L)$ is the supremum of the radii $R>0$ such that for any $v\in B_R(0)$ and for any $a\in H$, 
\be
\|av-v\|\leq L(a)
\ee
\end{enumerate}
\ed

\br\label{holradupbnd} In \eqref{propP}, by setting $u=0$ and $w=av$, it follows that the holonomy radius at zero is smaller than or equal to the convexity radius. Also, it follows that the convexity $\CvxRad$ radius satisfies that

\be
\CvxRad\leq\inf_{a\in H}\frac{L(a)}{\|a-id_V\|}
\ee
where $\|\cdot\|$ stands for the operator norm. It can be proved to be equal (see \cite{1004.1609}). This in turn gives an upper bound for the holonomy radius at zero; namely,
\be\label{HolRadineq}
\HolRad(0)\leq\inf_{a\in H}\frac{L(a)}{\|a-id_V\|}.
\ee
\er

\br\label{finflat} Notice that either of the radii given in Definition \ref{hlsprad} is finite if and only if $H\neq\{id_V\}$.
\er

In \cite{1004.1609}, the author proves that the restricted metric on the fibers of vector bundle are indeed holonomic spaces; namely, the following results.

\bt[\citet{1004.1609}]\label{lengthproof} Let $Hol_p$ be the holonomy group over a point $p\in M$ of a bundle with metric and connection and suppose that $M$ is Riemannian. Then the function $L_p:Hol_p\rightarrow \R$, 

\be\label{lengthnorm}
L_p(A)=\inf\{\ell(\alpha)|\alpha\in\Omega_p,P^{\alpha}_1=A\},
\ee
is a group-norm for $Hol_p$. Here $\Omega_p$ denotes the space of piece-wise smooth loops based at $p$. 
\et
\bd[\citet{1004.1609}]\label{lengthdef}
The function $L_p$, defined by \eqref{lengthnorm} will be called {\em length norm} of the holonomy group induced by the Riemannian metric at $p$.
\ed

\bt[\citet{1004.1609}]\label{holfib} Let $E_p$  be the fiber of a vector bundle with metric and connection $E$ over a Riemannian manifold $M$ at a point $p$. Let $Hol_p$ denote the associated holonomy group at $p$ and let $L_p$ be the group-norm given by \eqref{lengthnorm}. Then $(E_p,Hol_p, L_p)$ is a holonomic space. Moreover, if $E$ is endowed with the corresponding Sasaki-type metric, the associated holonomic distance coincides with the restricted metric on $E_p$ from $E$.
\et

\br\label{Eucdistni} Because the induced metric on any fiber is totally geodesic and flat, the identity map on $E_p$ between the Euclidean distance and the restricted distance is distance non-increasing (cf. Theorem \ref{hlmetthm}). 
\er

\bd\label{holraddef}  Given a vector bundle $\pi:E\rightarrow M$ with metric connection $\nabla$ and bundle metric $h$ over a Riemannian manifold $(M,g)$, the {\em holonomy radius} of $p\in M$ with respect to $\nabla$ is the holonomy radius at zero of the holonomic space $(E_p,Hol_p(\nabla),L_p)$.
\ed

\br By construction, the holonomy radius is positive; furthermore, by Remark \ref{finflat}, it is finite unless $Hol_p(\nabla)=\{id_{E_p}\}$. Notice that this is stronger than the connection being flat (since that would only imply that $Hol_p(\nabla)$ is zero dimensional).
\er

\br Similarly, the holonomy radius can be intrinsically defined for a Riemannian manifold by considering the one associated to the Levi-Civita connection (cf. \cite{1004.1609}). 
\er

\section{Convergence of Holonomic spaces}\label{GHhol}
Because holonomic spaces arise as fibers of vector bundles with metric connections over Riemannian manifolds (by Theorem \ref{holfib}, cf. \citet{1004.1609}), the study of their convergence properties becomes natural when trying to understand the behavior of said bundles endowed with their metrics of Sasaki-type under limits.  Also, given the underlying linear nature of the holonomic spaces, a $C^0$-convergence of the metrics to semimetrics is obtained, which implies the pointed Gromov-Hausdorff convergence of the holonomic spaces to precisely described spaces. 
Metrically, the description of their induced limit metrics is slightly more elusive.

\bt\label{holspcpt} Given a finite dimensional vector space, the collection of all holonomic space metrics $(V,d_L)$  is precompact in the $C^0$ sense. Namely, for any sequence $(V,H_i,L_i)$ there exists a subsequence (denoted without loss of generality with the same index $i$) for which the metrics $d_{L_i}:V\times V\rightarrow\R$ converge uniformly on bounded domains to a semi-metric $\rho:V\times V\rightarrow\R$.   
\et
\begin{proof}The strategy is the following: first use Arzel\`a-Ascoli on balls of a fixed radius $r>0$ around the origin; next argue that these convergences can be made to agree on $V$; and finally, argue that the limit function is a semi-metric. 

Let $V$ be a finite dimensional normed vector space. For any $r>0$ let $\varepsilon>0$ and consider $\eta=\min\{1,\frac{\varepsilon^2}{(1+2\sqrt{r})^2}\}$, $\delta=\frac{\eta}{\sqrt{2}}$.

Let $u,v, u',v' \in V$ be such that $\|u\|,\|v\|,\|u'\|,\|v'\|<r$ and
\be\label{delta}
\sqrt{\|u-u'\|^2+\|v-v'\|^2}<\delta.
\ee
Then, for any normed preserving linear map $a:V\rightarrow V$,
\begin{align}
\bigg| \|au-v\|- \|au'-v'\|\bigg | &\leq\|au-au'\|+\|v-v'\|\\ &\leq \|u-u'\|+\|v-v'\| \\ &<\sqrt{2}\sqrt{\|u-u'\|^2+\|v-v'\|^2}<\eta,\label{trineq}
\end{align}
by the triangle inequality for $\|\cdot\|$ and because $\|au-au'\|=\|u-u'\|$, $a\in O(V)$, \eqref{delta} and $\delta=\frac{\eta}{\sqrt{2}}$.

In particular, 
$$\|au-v\|^2\leq \|au'-v'\|^2+\eta^2+4r\eta,$$
which follows by direct squaring and by noticing that $\|au-v\|\leq 2r$, by the triangle inequality. 

Consider now any group-norm $L:H\rightarrow\R$ for any $H\leq O(V)$. By adding $L^2(a)$ to both sides, one sees that, 
$$L^2(a)+\|au-v\|^2\leq L^2(a)+\|au'-v'\|^2+\eta^2+4r\eta$$
now, by taking the square root and applying the triangle inequality, 

$$\sqrt{L^2(a)+\|au-v\|^2}\leq\sqrt{L^2(a)+\|au'-v'\|^2}+\eta+2\sqrt{r\eta}.$$

Now, since $\eta$ is less than one it follows that $\eta\leq\sqrt{\eta}$,  so that 
$$\sqrt{L^2(a)+\|au-v\|^2}\leq\sqrt{L^2(a)+\|au'-v'\|^2}+\varepsilon$$
holds. 

Therefore:
\begin{align}
d_L(u,v)&=\inf_a\sqrt{L^2(a)+\|au-v\|^2}\\ &\leq\inf_a\sqrt{L^2(a)+\|au'-v'\|^2}+\varepsilon\\ &= d_L(u',v')+\varepsilon.
\end{align}

Because of \eqref{trineq}, and by interchanging $u, v$ with $u',v'$, it now follows that
\be
|d_L(u,v)-d_L(u',v')|<\varepsilon,
\ee
thus proving that the family $\{d_L\}$ is equicontinuous on balls of a fixed radius $r>0$ around the origin in $V$. 

To prove uniform boundedness, one needs to observe that for any $L$ the following is true:
\be
d_L(u,v)\leq \sqrt{L^2(id_V)+\|id_Vu-v\|^2}=\|u-v\|\leq 2r.
\ee

The hypotheses of the classical  Arzel\`a-Ascoli's theorem now apply to get a uniform limit on the ball of radius $r>0$ (times itself). Consider now a countable exhaustion of $V$ by balls of radius $r_i\rightarrow\infty$. By a diagonal argument for any sequence of metrics $\{d_{L_i}\}$  one gets a pointwise limit $\rho$ on $V$ that is uniform on compact sets. 

Finally, except for nondegeneracy, all the properties of (semi)metrics are well behaved under limits:
\begin{align}
\rho(u,v)&=\lim_{i\rightarrow\infty}d_{L_i}(u,v)\geq0,\\ \notag\\
\rho(u,v)&=\lim_{i\rightarrow\infty}d_{L_i}(u,v)=\lim_{i\rightarrow\infty}d_{L_i}(v,u)=\rho(v,u),\\ \notag \\
\rho(u,u)&=\lim_{i\rightarrow\infty}d_{L_i}(u,u)=0,\\ \notag\\
\rho(u,v)&=\lim_{i\rightarrow\infty}d_{L_i}(u,v)\\Ê&\leq\lim_{i\rightarrow\infty}d_{L_i}(u,w)+\lim_{i\rightarrow\infty}d_{L_i}(w,v)\\ &=\rho(u,w)+\rho(w,v).
\end{align}

Therefore for any family of holonomic spaces $\{(V,d_{L_i)}\}$ there exists a subsequence that converges uniformly on compact sets.
\end{proof}

\br Nowhere in the proof was the fact that the function $d_L:V\times V\rightarrow\R$ was nondegenerate used; only properties of semi-metrics were required. However, the restriction to holonomic spaces yields nondegeneracy of $d_L$ and is of interest for the sequel. In general, the limit $\rho$ will be degenerate unless further assumptions are made (see Theorem \ref{holradbnd}).
\er

\bc\label{innholspcpt}Given a finite dimensional normed vector space, the collection of all pointed holonomic space metrics $((V,d_L), 0)$  is precompact in the pointed Gromov-Hausdorff sense. 
\ec
\begin{proof} By Theorem \ref{holspcpt} for any sequence $(V,d_i)$ of holonomic space metrics there exists a subsequence for which the semi-metrics $d_i$ converge uniformly on compact sets to a semi-metric $\rho$ on $V$.  The quotient space $Q=V/\sim$, where $u\sim v$ if and only if $\rho(u,v)=0$ for $u,v\in V$, is naturally a metric space; the metric is given by the distance $d_{\infty}([u],[v])=\rho(u,v)$ for any choice of representatives (or as the usual ---not Hausdorff--- distances between subsets of $V$). Therefore the convergent subsequence of metrics yields a convergent sequence of  metric spaces $(V,d_i)\rightarrow (Q,d_{\infty})$.
\end{proof}

\bc\label{precptHSGH} The space of holonomic metrics on inner-product spaces of dimension at most $k$ is precompact in the Gromov-Hausdorff sense. More explicitly, for any family of holonomic spaces $\{(V_i,H_i,L_i)\}$, where $V_i$ has dimension at most $k$ and its norm is induced by an inner product, there exists a subsequence that converges in the pointed Gromov-Hausdorff sense.
\ec
\begin{proof} By passing to a subsequence, one can assume that all the vector spaces in the sequence have the same dimension. Now, by Sylvester's Law of Inertia ---which in particular states that any two positive definite symmetric bilinear forms on a finite dimensional vector space are isometric---,  all the norms can be made to coincide, by way of some isometries $\phi_i:V\rightarrow V$. By defining $\tilde H_i=\phi_i\inv H_i\phi_i$ and $\tilde L_i:\tilde H_i\rightarrow\R$ by $\tilde L_i(b)=L_i(\phi_ib\phi_i\inv)$, the sequence $\{(V,\tilde H_i, \tilde L_i)\}$ now satisfies the hypotheses of  Theorem \ref{holspcpt}. 
\end{proof}

More can be said about the limiting metric spaces when there is more information about the underlying subgroups of isometries. Before that, recall the following result.

\bl[see \cite{MR1393940}. p. 47]\label{techlem} Let $(X,d)$ be a locally compact connected metric space and let $\{\varphi_i\}$ be a sequence of isometries of $(X,d)$. If there exists a point $x\in X$ such that $\{\varphi_i(x)\}$ converges, then there exists a subsequence $\{\varphi_{i_k}\}$ that converges to an isometry of $(X,d)$.
\el

\bt\label{limsup} Let $V$ be a finite dimensional normed vector space and $H$ be a subgroup of the group of norm preseving linear maps, denoted here by $O(V)$. Consider a sequence of group-norms $\{L_i:H\rightarrow\R\}$ such that the semi-metrics $d_i=d_{L_i}$, given by
\be\label{defholsem}
d_{L_i}(u,v)=\inf_a\sqrt{L_i^2(a)+\|au-v\|^2},
\ee
for any $u, v \in V$, converge uniformly on compact sets to a semi-metric $d_{\infty}$ on $V$. Then, there exists a set $G_0\subseteq O(V)$, given by: 
\be
G_0=\{g\in O(V)| g=\lim_{i\rightarrow\infty}a_i, \liminf_{i\rightarrow\infty}L_i(a_i)=0\},
\ee

such that for any $u,v\in V$,  
\be 
d_{\infty}(u,v)=0
\ee
if and only if there exists $g\in G_0$ such that $v=gu$, where $G_0$
\et
\begin{proof}
Let $u, v\in V$ be such that $d_{\infty}(u,v)=0$. This means that for any choice of $\varepsilon>0$ there exists $N=N_{\varepsilon}>0$ such that for any $j>N$,
\be
d_j(u,v)< \varepsilon
\ee

In particular by \eqref{defholsem} there exists $a_j(\varepsilon)\in H$ with 
\be
\sqrt{L_j^2(a_j)+\|a_ju-v\|^2}<\varepsilon,
\ee
which in turn gives that
\be
L_j(a_j),\|a_ju-v\|\leq\varepsilon.
\ee

By letting $\varepsilon=\frac{1}{n}$ and recursively choosing  $j=j_n=\max\{[N_{\frac{1}{n}}],j_{n-1}\}+1$, one produces a sequence $\{b_n=a_j(\frac{1}{n})\}$ for which $b_nu\rightarrow v$ and by Lemma \ref{techlem}, passing to a further subsequence if needed, such that it converges in $O(V)$ to some $g$, with $gu=v$ and $\lim_nL_{j(n)}(b_n)=0$.  The issue now is that this sequence is not parametrized by $i$ but by $n$.

To produce such a sequence, extend $\{b_n\}$ by
\be
c_i=b_n
\ee
for $i$ in the following range $j_n\leq i< j_{n+1}$. For this new sequence, the convergence 
$$c_i\rightarrow g$$ 
remains since only repeated terms were added in between the $b_n$'s, which were chosen to converge. Lastly, since for nonnegative sequences having vanishing limit inferior is equivalent to having a convergent subsequence converging to zero,  
$$\liminf_{i\rightarrow\infty}L_i(c_i)=0,$$
as promised, since letting $i_n=j_n$ provides the required subsequence.

Conversely, let $u\in V$ and consider $v=gu$ with $g\in G_0$, with $a_i\rightarrow g$.  Then for all $\varepsilon>0$ there exists $i>>0$ such that 
\be
d_i(u,gu)\leq \sqrt{L_i^2a_i+\|a_iu-gu\|^2}\leq \varepsilon.
\ee

So the claim now follows by the uniform convergence of $d_i\rightarrow d_{\infty}$. 
\end{proof}

\br It is not clear at this point whether $G_0$ is a subgroup of $O(V)$, since for example in the case of sequences of numbers $\{x_i\},\{y_i\}\subseteq\R$,
\be
\liminf_{i\rightarrow\infty}(x_i+y_i)\geq\liminf_{i\rightarrow\infty}x_i+\liminf_{i\rightarrow\infty}y_i.
\ee
Nevertheless, the characterization given by the previous theorem is still quite good as will be seen in the sequel. If one however insists upon having a group action to determine the degeneracy of $d_{\infty}$, this can be achieved by the following result. The drawback is that this new presentation says nothing about how to explicitly construct said group directly from the knowledge of $L_i$. 
\er

\bt\label{fbrqt} Let $V$ be a finite dimensional normed vector space and $H$ be a subgroup of the group of linear norm preserving isomorphisms, $O(V)$. Consider a sequence of group-norms $\{L_i:H\rightarrow\R\}$ such that the semi-metrics $d_i=d_{L_i}$, given by
\be\label{defholsemeq}
d_{L_i}(u,v)=\inf_a\sqrt{L_i^2(a)+\|au-v\|^2},
\ee
for any $u, v\in V$, converge uniformly on compact sets to a semi-metric $d_{\infty}$ on $V$. Then, there exists a closed subgroup of $O(V)$, given by
\be
G=\{g\in O(V)| \forall u\in V, d_{\infty}(u,gu)=0\},
\ee
such that for any $u,v\in V$ 
\be d_{\infty}(u,v)=0
\ee
if and only if there exists $g\in G$ such that $v=gu$.
\et
\br\label{G0CG}  Consider $g\in G_0$, as in Theorem \ref{limsup}. Then, for any $u\in V$, $d_{\infty}(u,gu)=0$ by Theorem \ref{limsup}. Thus 
$$G_0\subseteq G.$$ 
\er

\bd  The group $G$ will henceforth called the {\em wane group} of a convergent sequences of holonomic spaces. 
\ed

\begin{proof}[Proof of Theorem \ref{fbrqt}] Three statements must be proved: 1) the equivalence between having zero distance and being related by an element in $G$; 2) the fact that $G$ is actually a group; and 3) that this group is closed in $O(V)$.

To prove the equivalence first consider let $v\in V$ with $d_{\infty}(u,v)=0$, then by Theorem \ref{fbrqt} there exists $g\in G_0\subseteq G$  (by Remark \ref{G0CG}) with $v=gu$. 

Conversely, for any $g\in G$  and for any $u\in V$ 
\be
d_{\infty}(u,gu)=0,
\ee
by the definition of $G$. 

As the reader might have noticed, this doesn't prove that $G\subseteq G_0$ since this only implies that for any $u\in V$ and for any $g\in G$ there exists $h\in G_0$ such that 
\be
gu=hu.
\ee
Bear in mind that this equality is attained only at $u\in V$, since in principle $h$ depends on $u$. What was accomplished was the following: For any $u\in V$, $\{hu|h\in G_0\}=\{gu|g\in G\}$, that is that the equivalence classes determined by $G$ and by $G_0$ (in turn determined by the degeneracy of $d_{\infty}$) in $V$ are the same, as promised.

Secondly, to prove that $G$ is indeed a group, notice that because $d_{\infty}$ is already known to be a semi-metric,
\be
d_{\infty}(u,u)=0,
\ee
regardless of $u\in V$. So $id_V\in G$. 

Let $g,h\in G$, then by the triangle inequality of $d_{\infty}$,

\be
d_{\infty}(u,ghu)\leq d_{\infty}(u,hu)+d_{\infty}(hu,g(hu))=0+0
\ee
and, of course,
\be
d_{\infty}(u,g\inv u)=d_{\infty}(g(g\inv u),g\inv u)=0.
\ee
Therefore $G$ is a subgroup of the group of norm preserving linear maps, $O(V)$.

Finally, to see that $G$ is closed, notice that for any $u\in V$ the assignment $\varphi_u:O(V)\rightarrow\R$, given by 
\be
\varphi_u:g\mapsto d_{\infty}(u,gu),
\ee
is a composition of continuous functions and thus itself continuous. Because of this, $G$ can be represented as the following intersection of closed sets,
\be
G=\bigcap_{u\in V}\varphi_u\inv(0),
\ee
and thus $G$ is closed.
\end{proof}
\br  If $V$ is further assumed to be an inner product space, then $O(V)$ is a compact Lie group; furthermore, because by Theorem \ref{fbrqt}, $G$ is also a compact Lie group. 
\er

As mentioned before, the assumption that the sequence of holonomic space metrics have a single group $H$ can be weakened.

\bd\label{finreptyp} A family of subgroups $\{H_i\}$ of $O(V)$ is of finite representation type if the set  $\{H_i\}$  be finite up to isomorphism, and that this isomorphism can be chosen to be by conjugation by an element in $O(V)$. 
\ed

\bt\label{dinfty} Let $V$ be a normed vector spaces and consider a sequence of holonomic space metrics $\{d_i\}$ on it such that the collection of corresponding groups $\{H_i\}$ is of finite representation type. Suppose further that it is a convergent sequence. Then there exists a closed subgroup $G\leq O(V)$ (as in Theorem \ref{fbrqt}) such that for any $u,v\in V$ 
\be d_{\infty}(u,v)=0
\ee
if and only if there exists $g\in G$ such that $v=gu$.
\et
\begin{proof} The only subtlety to consider here is the fact that the groups need not be the same.  However, by assumption there exists a finite set $\{K_1,\dots,K_n\}$ of subgroups of $O(V)$, $\phi_i:V\rightarrow V$ linear isometries, and a function $\psi:\N\rightarrow \{1,\dots,n\}$ such that
\be H_i=\phi_i\inv K_{\psi(i)}\phi_i.
\ee

Define $\tilde L_i:K_{\psi(i)}\rightarrow\R$ as follows. For any $k\in K_{\psi(i)}$,
\be
\tilde L_i(k)=L_i(\phi_ik\phi_i\inv).
\ee

It remains to show that $(V,H_i,L_i)$ and $(V,K_{\psi(i)},\tilde L_i)$ yield isometric holonomic spaces. For this observe that
\begin{align}
d_{L_i}(u,v)&=\inf_{a\in H_i}\sqrt{L_i^2(a)+\|au-v\|^2}\\
&=\inf_{b\in K_{\psi(i)}} \sqrt{\tilde L_i^2(b)+\|\phi_i\inv b\phi_iu-v\|^2}\\
&=\inf_{b\in K_{\psi(i)}}\sqrt{\tilde L_i^2(b)+\|b\phi_iu-\phi_iv\|^2}\\
&=d_{\tilde L_i}(\phi_iu,\phi_iv).
\end{align}

Now, since there are only finitely many different groups, one can pass to a subsequence that consists of a single one. This is now the setting for Theorem \ref{fbrqt}. 
\end{proof}

\bc\label{inndinfty} Let $\{V_i\}$ be a collection of finite dimensional inner product vector spaces and consider a sequence of holonomic space metrics $\{d_i\}$ on $\{V_i\}$ such that the collection of corresponding groups $\{H_i\}$ consist of finite representation types. Suppose further that the sequence of metric spaces $\{(V_i,d_i)\}$ is a convergent sequence in the Gromov-Hausdorff sense. Then there exists a positive integer $k$ and a closed subgroup $G\leq O(k)$ such that 
\be 
(V_i,d_i)\ptGH\R^k/G,
\ee
where the metric on the limit is obtained as in Corollary \ref{precptHSGH}.
\ec
\begin{proof} As in Corollary \ref{precptHSGH}, one can pass to a subsequence and assume that $\{V_i\}$ has constant dimension $k$ and such that the norms are the constant. The conclusion now follows from Theorem \ref{dinfty}.
\end{proof}

Recall that (in Definition \ref{hlsprad}) for a holonomic space $(V,H,L)$ the {\em holonomic radius} at a point $u\in V$, $\HolRad(u)$, is the largest $r>0$ for which the metric $d_L$ is isometric to the Euclidean metric when restricted to the ball of radius $r$ around $u\in V$. In the special case when $u=0$, then the following result holds. 

\bl[see Remark \ref{holradupbnd}, cf. \citet{1004.1609}] Given a holonomic space $(V,H,L)$, 
\be
\HolRad(0)\leq\inf_{a\in H}\frac{L(a)}{\|a-id_V\|}.
\ee
\el

\bt\label{holradbnd}
Let $(V,H,L_i)$ be a convergent sequence of holonomic spaces. Suppose further that there exists a constant $c>0$ such that 
\be
c\leq\HolRad_i(0).
\ee
Then the limit semi-metric $d_{\infty}$ is nondegenerate. That is that $d_{\infty}$ is a metric.
\et
\begin{proof} Consider $g\in G_0$, as in \ref{limsup}, and let $\{a_i\}$ be any defining sequence for $g$. That is such that $a_i\rightarrow g$ and $\liminf_{i\rightarrow\infty}L_i(a_i)=0$. By Lemma \ref{holradupbnd}, for each $a_i$, 
\be
c\leq\HolRad_i(0)\leq\frac{L_i(a_i)}{\|a_i-id_V\|}.
\ee

Thus for any $\varepsilon>0$ and for any $N>0$ there exists $i>N$, 
\be
\|a_i-id_V\|\leq\frac{L_i(a_i)}{c}\leq\frac{\varepsilon}{c}.
\ee

Therefore there is a subsequence of $\{a_i\}$ that converges to the identity map. Because $\{a_i\}$ converges to $g$, it follows that $g=id_V$ and now Theorem \ref{limsup} yields the claim.
\end{proof}

Finally, as an application of these concepts the upper semi-continuity of the holonomy radius of a connection over a Riemannian manifold can be asserted (recall Definition \ref{holraddef}). In the case of a holonomic space $(V,H,L)$, the author proved in \cite{1004.1609} that the holonomy radius is a continuous function on $V$. 

\bp Given a vector bundle $\pi:E\rightarrow M$ with metric connection $\nabla$ and bundle metric $h$ over a Riemannian manifold $(M,g)$, then the holonomy radius $\HolRad:M\rightarrow \R$ is an upper semicontinuous function.
\ep
\begin{proof} Let $p\in M$ and consider a sequence $\{p_i\}\subseteq M$ converging to $p$. By Remark \ref{ptsqnHaus} and by the fact that the fibers of $\pi$ are equidistant, it follows that the holonomic metrics $\{d_{ L_{p_i}}\}$ converge in the $C^0$ sense to the holonomic metric $d_{L_p}$. 

Let $\rho=\HolRad(p)$. Now,  since the metrics $\{d_{ L_{p_i}}\}$ converge uniformly when restricted to the ball of radius $\rho$, and the metric $d_{L_p}$ is Euclidean on that ball, let 
$$\tilde\rho=\limsup_{i\rightarrow\infty}\HolRad(p_i).$$ 
Let $u,v\in E_p$ with $\|u\|,\|v\|<\tilde\rho$. Then, there exist (sub)sequences $\{u_i\}\subseteq E_{p_i}$, $\{v_i\}\subseteq E_{p_i}$ with $\|u_i\|,\|v_i\|<\tilde\rho$, converging to $u$ and $v$ respectively such that 
\be
d_{L_{p_i}}(u_i,v_i)=\|u_i-v_i\|.
\ee
Therefore,
\be
d_{L_p}(u,v)=\lim_{i\rightarrow\infty}d_{L_{p_i}}(u_i,v_i)=\lim_{i\rightarrow\infty}\|u_i-v_i\|=\|u-v\|,
\ee
thus proving that $\tilde\rho\leq\rho$ as promised.
\end{proof}

\section{Gromov-Hausdorff convergence of metrics of Sasaki-type}\label{GH}

Since vector bundles (with metric connections) appear naturally as associated objects to Riemannian manifolds, providing the latter with additional structures (as Poisson brackets, control structures, orientations, holomorphic structures, etc.), it is natural to investigate the behavior of these bundles under limits of their bases. 

One reason to analyze them from the viewpoint of holonomic spaces is given by Theorem \ref{holfib} (see also \citet{1004.1609}). Furthermore, the next result gives yet another reason.

\bt\label{incGH} Given a vector bundle $\pi:E\rightarrow M$ with metric connection $\nabla$ and bundle metric $h$ over a Riemannian manifold $(M,g)$, consider a point $p\in M$ and let $(V,H, L)$ be the holonomic space $(E_p,Hol_p(\nabla),L_p)$. Then the Gromov-Hausdorff distance between $(V,d_L)$ and $\pi\inv(B_R(p))\subseteq E$ (with the restricted metric from $E$) is finite and bounded by $2R$.
\et
\begin{proof}
By \cite{1004.1609}, the inclusion $(E_q,d_{L_q})\hookrightarrow E$ is an isometric embedding in the sense of metric spaces for any q. Furthermore, because the projection map is a Riemannian submersion, parallel translation along any minimal geodesic in $M$ connecting the points $p,q\in M$ renders the fibers equidistant.  Therefore, the distance between the central fiber and any other fiber over a ball of radius $R$ is bounded by $R$.  From this, for any $p\in M$ and $R>0$ the inclusion map of the central fiber 
\be
E_p\hookrightarrow \pi\inv(B_R(p))\subseteq E
\ee
is an $R$-isometry. By Proposition \ref{epsiso}, the claim now follows. 
\end{proof}

In particular, for the tangent bundle: 

\bc Given a Riemannian manifold $(M,g)$ and a point $p\in M$ let $(V,d_L)$ be the holonomic space $(M_p,Hol_p(g),L_p)$ then the Gromov-Hausdorff distance between $(V,d_L)$ and $\pi_M\inv(B_R(p))\subseteq TM$ is bounded by $2R$.
\ec

\bt \label{precptBWC} Given a precompact collection of (pointed) Riemannian manifolds $\mathcal{M}$ and a positive integer $k$, the collection $\textup{BWC}_k(\mathcal{M})$ of vector bundles with metric connections of rank $\leq k$ endowed with metrics of Sasaki-type is also precompact. The distinguished point for each such bundle is the zero section over the distinguished point of their base.
\et
\br Passing to a subsequence is unavoidable as can be seen in Example \ref{prodcollapse}.
\er
\begin{proof}[Proof of Theorem \ref{precptBWC}] Fix $\varepsilon>0$ and $R>0$. Following Theorem \ref{Gromovs} and Remark \ref{GromovN}, define $C=C(\epsilon,R)>0$ by 
\be 
C:=\max_{i\leq k}N(\varepsilon, R, \E^i),
\ee  
where $\E^i$ is the Euclidean space of dimension $i$. 

Let $(E,h,\nabla)\stackrel{\pi_E}{\longrightarrow} (M,g)$ be any bundle with metric connection and let $N(R,\varepsilon)$ be the uniform bound on the number of balls of radius $\varepsilon$ needed to cover a ball of radius $R$ on $\mathcal{M}$. Consider  $p\in M$ and its zero section $0=\varsigma(p)$. 

Since $\pi_E$ is a Riemannian submersion, $\pi_E(B_R(0))=B_R(p)$. Let $A$ be any $\varepsilon$-net in $B_R(p)$.  Since for each $a\in A$, $E_p$ is flat, let $A_a$ be any $\varepsilon$-net in $B_R(\varsigma(a))\subseteq T_aM$. The cardinality of $A_a$ can be chosen to less than $C(\varepsilon,R)$, because the identity map is a distance non-increasing map between the (induced) Euclidean metric on $E_p$ and the restricted metric (see Remark \ref{Eucdistni}).  

Let $u\in B_R(0)$, then let $a\in A$ such that $d(a,\pi_E(u))<\varepsilon$. Let $\gamma$ be any minimal geodesic connecting $\pi_mu$ to $a$ and let $v\in T_aM$ be the parallel image of $u$ along $\gamma$. Finally consider $u_a\in A_a$ to be such that $|u_a-v|<\varepsilon$. Hence,
\[
d(u,u_a)\leq d(\pi_Mu,a)+d(v,u_a)\leq2\varepsilon.
\]
Therefore, given any $R>0$ , and for any $\varepsilon>0$, the following holds.

\be
N(2\varepsilon, R, E)\leq N(\varepsilon, R, M)+ C(\varepsilon,R).
\ee
 
So that if the assignment $M\mapsto N(\varepsilon, R, M)$ is bounded on $\mathcal{M}$, then so is $E\mapsto N(\varepsilon, R, E)$ on $\textup{BWC}_k(\mathcal{M})$.
 
Therefore, in view of Gromov's Compactness Theorem \ref{Gromovs}, this finishes the proof.
\end{proof}


\bp\label{sascpt} For any sequence of Riemannian manifolds $\{(X_i,p_i)\}$ converging to $(X_{\infty},x_{\infty})$ consider a convergent family of bundles with metric connection $(E_i,h_i,\nabla_i)$ over it converging to $(E_{\infty},y_{\infty})$. Then there exist continuous maps $\pi_{\infty}:E_{\infty}\rightarrow X_{\infty}$, $\varsigma_{\infty}:X_{\infty}\rightarrow E_{\infty}$, $\mu_{\infty}:E_{\infty}\rightarrow\R$, and a subsequence, without loss of generality also indexed by $i$, such that:
\begin{enumerate}
\item the projection maps $\pi_i:E_i\rightarrow X_i$ converge to $\pi_{\infty}:E_{\infty}\rightarrow X_{\infty}$, which is also a submetry;
\item the zero section maps $\varsigma_i:X_i\rightarrow E_i$ converge to $\varsigma_{\infty}:X_{\infty}\rightarrow E_{\infty}$, which is also a isometric embedding; 
\item $\pi_{\infty}\circ \varsigma_{\infty}=id_{X_{\infty}}$; and
\item the maps $\mu_i:E_i\rightarrow\R_{\geq0}$, given by 
$$\mu_i(u)=d_{E_i}(u,\varsigma_i\circ\pi_i(u))=\sqrt{h_i(u,u)},$$ 

converge to $\mu_{\infty}:E_{\infty}\rightarrow\R_{\geq0}$, also given by $$\mu_{\infty}(y)=d_{E_{\infty}}(y,\varsigma_{\infty}\circ\pi_{\infty}(y)).$$
\end{enumerate}
\ep

\begin{proof} Since all of these maps preserve the distinguished points (consider $0\in\R_{\geq0}$), by the Arzel\`a-Ascoli Theorem \ref{ArAs}, one only has to check equicontinuity. But this is immediate from the fact that the $\pi_i$ are submetries \cite{MR2500106}, $\varsigma_i$ isometric embeddings, and $\mu_i$ both distance functions and submetries.  In fact, their limits will share these properties, as noted by \citet[Section 10.1.3]{MR2243772}.

The equation $\pi_{\infty}\circ \varsigma_{\infty}=id_{X_{\infty}}$ holds since the corresponding equation holds for every $i$.  Finally, the equation $\mu_{\infty}(y)=d_{E_{\infty}}(y,\varsigma_{\infty}\circ\pi_{\infty}(y))$ holds, since for any sequence $\{u_i\}$ converging to $y\in Y$ the geodesics $t\mapsto tu_i=\varsigma_i\circ\pi_i(u_i)+tu_i$ are rays (see Remark \ref{sasprops}), hence isometric embeddings, and thus also converge to a minimal geodesic. 
\end{proof}

Because of Theorem \ref{incGH}, the fiberwise behavior is also controlled. 

\bp\label{fibsGH} Let $\pi_i:(E_i,0_{*_i})\rightarrow (X_i,*_i)$ be a convergent sequence of pointed spaces as before. Let $\pi_{\infty}:(E_{\infty},0_{*_{\infty}})\rightarrow (X_{\infty},*_{\infty})$ be their limit. Then if $q\in X_{\infty}$ and $\{q_i\in X_i\}$ is any sequence converging to $q$. Then, by passing to a subsequence if needed, for any $\varepsilon>0$, 
\be\label{epsconv} 
\pi_i\inv(B_{\varepsilon}(q_i))\ptGH\pi_{\infty}\inv(B_{\varepsilon}(q)).
\ee
Furthermore, 
\be\label{fibconv} 
\pi_i\inv(q_i)\ptGH\pi_{\infty}\inv(q).
\ee
\ep
\begin{proof}Since the pointed sequence converges with distinguished point $*_i$, it also converges with respect to the points $q_i$.  

Consider, as in the proof of Proposition \ref{sascpt}, the minimizing geodesics $\gamma_i$ given by  $t\mapsto tu_i$  for any convergent sequence of points $u_i\in\pi_i\inv(q_i)$. Then, the sequence converges to a minimizing geodesic and since the sequence of maps $\pi_i\circ\gamma_i\equiv q_i$ also converges, it follows that the limit $Q$ of the fibers (which is is known to exist by Theorem \ref{holspcpt} or Corollary \ref{innholspcpt} ) is inside the fiber over the limit (see Remark \ref{ptsqnHaus}). More precisely, there exists an isometric embedding 

\be Q\hookrightarrow \pi_{\infty}\inv(q).
\ee 

To prove that this is indeed surjective, and thus proving \eqref{fibconv}, the statement of \eqref{epsconv} will be proved first.

For any $\varepsilon>0$, the sequence $(\pi_i\inv(B_{\varepsilon}(0_{q_i})),0_{q_i})$ also converges (or a subsequence thereof) by Theorem \ref{precptBWC}.  Any sequence of points $u_i\in\pi_i\inv(B_{\varepsilon}(q_i))$ that converges, necessarily converges to a point $y\in \pi_{\infty}\inv(B_{\varepsilon}(q))$ since there exists $N>0$ such that for any $i>N$ 
\be
d_i(q_i, \pi_i(u_i))\leq\varepsilon,
\ee
so that, by continuity, 
\be
d_{\infty}(q, \pi_{\infty}(y))=\lim_{i\rightarrow\infty}d_i(q_i, \pi_i(u_i))\leq\varepsilon
\ee

Conversely, consider any $y\in \pi_{\infty}\inv(B_{\varepsilon}(q))$ and any sequence $\{u_i\}$ converging to $y$. By looking again at $\gamma_i$, the minimizing geodesics from $\varsigma_i\circ\pi_i(u_i)$ to $u_i$, one sees that a subsequence of $\{\varsigma_i\circ\pi_i(u_i)\}$ converges to $\varsigma_{\infty}\circ\pi_{\infty}(y)$. By Proposition \ref{sascpt}, $\pi_{\infty}$ is an isometry when restricted to the image of $\varsigma_{\infty}$; therefore there exists a subsequence of $\{\pi_i(u_i)\}$ that converges to $\pi_{\infty}(y)$. Now, because 
\be
\pi_{\infty}(y)\in B_{\varepsilon}(q),
\ee
it follows that there exists $N>0$ such that for all $i>N$,
\be
\pi_i(u_i)\in B_{\varepsilon}(q_i).
\ee

Thus, for all $\varepsilon>0$, 

\be
\pi_i\inv(B_{\varepsilon}(q_i))\ptGH\pi_{\infty}\inv(B_{\varepsilon}(q)).
\ee

Furthermore,  this convergence is attained in a compatible way with the convergence of their ambient spaces.

To finish the proof of \eqref{fibconv} consider any $y\in\pi_{\infty}\inv(q)$ and any sequence $\{u_i\}$ converging to $y$.  Since
$$
\pi_{\infty}\inv(q)\subseteq \pi_{\infty}\inv(B_{\varepsilon}(q))
$$
for any $\varepsilon>0$, by \eqref{epsconv} the sequence can be assumed to satisfied that 
$$u_i\in\pi_i\inv(B_{\varepsilon}(q_i)).$$
It is better to denote this sequence by $\{u^{\varepsilon}_i\}$, since it in fact depends on $\varepsilon$. Now, by Theorem \ref{incGH}, for any $\varepsilon>0$,
$$d_{GH}(\pi_i\inv(B_{\varepsilon}(0_{q_i})), \pi_i\inv((q_i)))<2\varepsilon.$$
One can consider a sequence $\{\widetilde{u^{\varepsilon}_i}\in\pi_i\inv(q_i)\}$ with
\be
d_{E_i}(u^{\varepsilon}_i,\widetilde{u^{\varepsilon}_i})<\varepsilon.
\ee

Now, again by a diagonalization argument, consider $\varepsilon=\frac{1}{i}$ and define

\be v_i=\widetilde{u^{\varepsilon}_i}.
\ee

By definition, $v_i\in\pi_i\inv(q_i)$ and for any $\varepsilon>0$, there exist $N$ for such that for any $i>N$, $$2/i<\varepsilon,$$
and as such, 

\be d_{E_i}(u^{1/i}_i,v_i)<\frac{\varepsilon}{2},
\ee 
and
\be
d(u^{1/i}_i,y)<\frac{\varepsilon}{2}.
\ee
Therefore, for any $y$ over $q$, a sequence  $\{v_i\}$ over $q_i$ that converges to $y$ was produced. By Remark \ref{ptsqnHaus}, $Q\cong \pi_{\infty}\inv(q)$ and the claim follows.
\end{proof}

\bp \label{equidfib} For any sequence of Riemannian manifolds $\{(X_i,p_i)\}$ converging to $(X_{\infty},x_{\infty})$ consider a convergent family of bundles with metric connection $(E_i,h_i,\nabla_i)$ over it converging to $(E_{\infty},y_{\infty})$  and $\pi_{\infty}:E_{\infty}\rightarrow X_{\infty}$ as in Proposition \ref{sascpt}. Then the fibers of $\pi_{\infty}$ are equidistant.
\ep
\begin{proof}
Let $p,q\in X_{\infty}$ be arbitrary and consider sequences $\{p_i\},\{q_i\}\subseteq X_i$ converging to $p,q\in X_{\infty}$ respectively. Because for each $i$ the fibers of $\pi_i$ are equidistant,  the distance between the fibers $\pi_i\inv(p_i)$ and $\pi_i\inv(q_i)$ is equal to the distance between $p_i$ and $q_i$. 

Let $u,v\in E_{\infty}$ and consider sequences $\{u_i\in\pi_i\inv(p_i)\}$ and $\{v_i\in\pi_i\inv(q_i)\}$ converging to $u,v$ respectively. Then,

\be
d_{E_i}(u_i,v_i)\geq d_{X_i}(p_i,q_i);
\ee
which in turn implies that

\be
d_{E_{\infty}}(u,v)\geq d_{X_{\infty}}(p,q).
\ee

This proves that the distance between the fibers is at least the distance between their base points. 

It remains to show that for any $u\in\pi_{\infty}\inv(p)$ there exists $v\in\pi_{\infty}\inv(q)$ with  

\be
d_{E_{\infty}}(u,v)=d_{X_{\infty}}(p,q).
\ee

To see this, consider any sequence $\{u_i\in\pi_i\inv(p_i)\}$ converging to $u$. Let $\{v_i\in\pi_i\inv(q_i)\}$ be a sequence such that
\be
d_{E_i}(u_i,v_i)= d_{X_i}(p_i,q_i).
\ee

Let $\alpha_i$ be minimizing geodesics connecting $u_i$ to $v_i$. By the Arzel\`a-Ascoli Theorem, there exists a subsequence of $\{\alpha_i\}$ that converges to a minimizing geodesic in $E_{\infty}$ connecting $u$ to some point $v\in E_{\infty}$. It follows that the corresponding subsequence of $\{v_i\}$ converges to $v$. From this it follows, since fibers converge to fibers, that for any $p,q\in X_{\infty}$ and for any $u\in\pi_{\infty}\inv(p)$ there exists 
\be
v\in\pi_{\infty}\inv(q),
\ee 
such that, by continuity, 

\be
d_{E_{\infty}}(u,v)=\lim_{i\rightarrow\infty}d_{E_i}(u_i,v_i)=\lim_{i\rightarrow\infty} d_{X_i}(p_i,q_i)=d_{X_{\infty}}(p,q).
\ee
\end{proof}

Suppose further that $\{E_i\}$ is of finite holonomy representation type (see Definition \ref{finreptyp}), then the fibers of $\pi_{\infty}$ can be naturally identified with the quotient of  any given fiber by a closed subgroup of the orthogonal group in view of Theorem \ref{inndinfty}, as stated in the following result.

\bt\label{finholtype} Let $\pi_i:E_i\rightarrow X_i$ be a convergent sequence of vector bundles with bundle metric and compatible connections $\{(E_i,h_i, \nabla_i)\}$, with limit $\pi:E\rightarrow X$. Suppose further that there are only finitely many holonomy representation types. Then there exists a positive integer $k$ such that for any point $p\in X$ there exists a compact Lie group $G\leq O(k)$, henceforth called the {\em wane group}, that depends on the point, such that the fiber $\pi\inv(p)$ is homeomorphic to $\R^k/G$, i.e. the orbit space under the standard action of $G$ on $\R^k$.
\et
\br Recall that the main feature of $G$ is that its orbits coincide with the ``orbits'' of the following set $G_0$ (see Theorems \ref{limsup} and \ref{dinfty}). Let $p\in X$ and let $p_\in X_i$ be a sequence that converges to $p$. Suppose, by passing to a subsequence, that the sequence of holonomy groups is constant, $H\equiv Hol_{p_i}(\nabla_i)$, and the fibers have constant dimension. Then let $G_0$ is given by
\be
G_0=\left\{g=\lim_{i\mapsto\infty}a_i\left|a_i\in H, \liminf_{i\mapsto\infty}L_i(a_i)=0\right.\right\},
\ee
where $L_i(a_i)$ is the infimum of the lengths of loops at $p_i\in X_i$ that generate $a_i$ by parallel translation (as in Definition \ref{lengthdef}).

By virtue of Corollaries \ref{precptHSGH} and \ref{inndinfty},  the limit metric of $\pi\inv(p)$ could, in principle, be given more explictly, once the behavior of these lengths is known. 
\er

\begin{proof}[Proof of Theorem \ref{finholtype}]  

Let $\{p_i\}$ be a sequence converging to $p$. Then, by Proposition \ref{fibsGH}, there exists a subsequence of $\{\pi_i\inv(p_i)\}$ that converges to the fiber $\pi\inv(p)$.  Let $V_i=\pi_i\inv(p_i)$, $H_i=Hol_{p_i}(\nabla_i)$, and $L_i:H_i\rightarrow\R$ the induced length norm. Then by Corollary \ref{inndinfty}, applied to the sequence $\{(V_i,H_i,L_i)\}$, the conclusion now follows.
\end{proof}

Summarizing, in the case of the collection of Sasaki metrics on the tangent bundles of a convergent sequence of simply connected Riemannian manifolds, the following holds.

\bt Let $\{M_i\}$ be a family of simply connected Riemannian manifolds that converges in the (pointed) Gromov-Hausdorff sense to $X$. Then there exists a subsequence of $\{TM_i\}$, with their Sasaki metrics, that converges to a space $Y$. Furthermore,  
\begin{enumerate}
\item there exists a continuous map $\pi:Y\rightarrow X$, that is a submetry with equidistant fibers.
\item there exists a positive integer $k$ such that for any $p\in X$ there exists a closed subgroup $G\leq O(k)$ such that $\pi\inv(p)$ is homeomorphic to $\R^k/G$.
\end{enumerate}
\et
\begin{proof} 
Because the sequence $\{M_i\}$ is convergent, by Theorem \ref{precptBWC}, there is a subsequence of $\{TM_i\}$ that converges to a space $Y$. By Propositions \ref{sascpt} and \ref{equidfib}, the promised $\pi:Y\rightarrow X$ exists and has the required properties Finally, because the manifolds are simply connected, by virtue of the classification theorem of \citet{MR0079806}), there exist finitely many holonomy representation types for the sequence. Therefore, by Theorem \ref{finholtype} the rest of the claim holds.
\end{proof}

Theorems \ref{holradbnd} and \ref{finholtype} together give a criterion for the fibers of $\pi_{\infty}$ in Theorem \ref{sascpt} to be vector spaces:

\bt\label{hlradbnd}  Let $\pi_i:E_i\rightarrow X_i$ be a convergent sequence of vector bundles with bundle metric and compatible connections $\{(E_i,h_i, \nabla_i)\}$, with limit $\pi:E\rightarrow X$. Suppose further that there are only finitely many holonomy representation types (cf. Theorem \ref{inndinfty}) and that there exist a uniform positive lower bound for the holonomy radii of $\pi_i:E_i\rightarrow X_i$ as in Definition \ref{holraddef}. Then the fibers of  $\pi_{\infty}$ are vector spaces. 
\et
\begin{proof} Again, by reduction to the case where the rank is constant and where there is a single holonomy representation, the conclusion follows from Theorem \ref{holradbnd}.
\end{proof}

Another consequence of Theorems \ref{sascpt} and \ref{fbrqt} is the following. 

\bx\label{totcollapse} Given any compact Riemannian manifold $(M^n,g)$ and let $X_i$ be the metric space obtained by rescaling $g$ into $\frac{1}{i^2}g$. Then tangent spaces converge to $\R^n/\overline{Hol(g)}$. Here $\overline{Hol(g)}$ denotes the closure in $O(n)$.
\ex
\begin{proof} Let $p\in M$ and let $(V,H,L)$ be the associated holonomy space given by Theorem \ref{holfib} at $p$.  By Theorem \ref{incGH}, because $M$ is compact, it follows that 
\be
d_{GH}(TM, V)\leq 2\diam(M).
\ee

Thus, by re-scaling, $\diam(X_i)\rightarrow0$ and the limit $Y=\lim_iTX_i$ is equal to the limit of the holonomic metrics at $p$. To analyze these spaces notice that the Sasaki metric re-scales like the base metric does; indeed, by re-scaling, the Levi-Civita connection remains constant and thus the horizontal lifts remain unchanged (cf. Definition \ref{sasmetdef}). Define $$(V_i,H_i,L_i)$$ to be the corresponding holonomic spaces at $p\in X_i$.  Again, because the connection is unchanged, it follows that $$H_i\equiv H.$$ 

Also, since by Definition \ref{lengthdef} $L$ is an infimum of lengths,  
$$L_i=\frac{1}{i}L.$$

Finally, the norm on $V_i$, denoted by $\|\cdot\|_i$, is given by 
$$\|\cdot\|_i=\frac{1}{i}\|\cdot\|,$$

where $\|\cdot\|$ is the norm on $V$.  Notice that by considering the map $\phi_i:V\rightarrow V$, given by
\be \phi_i: u\mapsto iu,
\ee
one gets an isometry of the holonomic spaces
\be
\phi_i:(V,H,L_i)\rightarrow (V_i,H_i, L_i).
\ee

Therefore, the limit $Y$ is the quotient $\R^n/G$ for some compact Lie group $G$, by Theorem \ref{fbrqt}. Furthermore, consider $G_0$ as in Theorem \ref{limsup}. Because for any $a\in H$, 
$$\lim_{i\rightarrow\infty}L_i(a)=0,$$
it follows that $G_0=\overline{H}$; thus proving the claim.  

\end{proof}
\br A theorem of \citet{MR1758296} states that any closed subgroup of $O(n)$ can be realized as the closure of a holonomy group of a compact smooth manifold. By Theorem \ref{totcollapse} one thus recovers all linear metric quotients of $\R^n$.
\er

\bx\label{prodex} Let $\{(M_i,g_i)\}$ and $\{(N_i,h_2)\}$ be two convergent sequences of complete Riemannian manifolds, with limits $X$ and $Y$ respectively. Let $\{E_i\rightarrow M_i\}$ and $\{F_i\rightarrow N_i\}$ be two convergent sequences of vector bundles with connections endowed with their metrics of Sasaki-type. Let $E\rightarrow X$ and $F\rightarrow Y$. Then, for the product metrics on $M_i\times N_i$ the limit converges to $E\times F\rightarrow X\times Y$.
\ex
\begin{proof} For each $i$ the bundles $E_i\times F_i\rightarrow M_i\times N_i$ are endowed with the product connection, the product bundle metric, from which it follows that for any curve $\gamma=(\gamma_1,\gamma_2)$, the parallel translation along $\gamma$ splits in the following way.
$$P^{\gamma}=P^{\gamma_1}\oplus P^{\gamma_2}$$

From this it follows that the product metric of metrics of Sasaki-type coincides with the metric of Sasaki-type on the product.  Now, because the spaces are products, the limit of the product is the product of the limits.
\end{proof}
\br
Because of this, it follows that for any $(p,q)\in M_i\times N_i$, 
\be
Hol_{(p,q)}(\tilde g_i)=Hol_p(E_i)\times Hol_q(F_i).
\ee
Furthermore, the length norm of $(a,b)$ is given by 
\be
L_{(p,q)}\left((a,b)\right)=\sqrt{L_p^2(a)+L_q^2(b)} 
\ee
\er
\bx\label{prodcollapse} Let $(M_1^{n_1},g_1)$ and $(M_2^{n_2},g_2)$ be two complete Riemannian manifolds; suppose further that $M_2$ compact. Consider the metrics $$\{\tilde g_i=g_1+\frac{1}{i^2}g_2\}$$ on $M_1\times M_2$, which to converge to $(M_1,g_1)$. If $\pi: Y\rightarrow M_1$ is the limit of their corresponding tangent bundles endowed with their metrics of Sasaki-type, then the fibers of $\pi$ are homeomorphic to 
\be
\R^{n_1}\times \left(\R^{n_2}/\overline{Hol(g_2)}\right).
\ee
However, for the constant sequence $\{(M_1,g_1)\}$ the limit is the canonical projection $TM_1\rightarrow M_1$. This proves that passing to a subsequence in Theorem \ref{precptBWC} is in general unavoidable.
\ex
\begin{proof}
Because these metrics are product metrics, where the second factor is re-scaled, the limit is the limit of the factors, by Example \ref{prodex}. 

Now, since only one of the factor is being re-scaled, while the other remains constant, the group-norm becomes degenerate on $\{id\}\times Hol(M_2)$, thus yielding the desired result by Example \ref{totcollapse}.
\end{proof}

\section{A weak notion of parallelism on singular spaces}
The most general setting for a notion of parallel translation is that of a pair of metric spaces and a surjective submetry between them. Even in this generality one can talk about parallel translation as long as one is willing to loosen it by considering, instead of functions, relations or ---equivalently--- set-valued functions. 

In the setting of limits of vector bundles with connection endowed with their metrics of Sasaki-type  the following general considerations will be considerably better behaved. Yet, in giving a precise framework the assumptions will be kept to a bare minimum.

\subsection{Parallelism for submetries} 

\bd Let $X$ and $Y$ be metric spaces a surjective {\em submetry} $\pi:Y\rightarrow X$ is a surjective map ({\em a fortiori} continuous) such that for any radius $r$  the image of any metric ball of radius $r$ is again a ball of the same radius $r$. 
\ed

\bd Given a surjective submetry $\pi:Y\rightarrow X$ a curve $\gamma:[0,1]\rightarrow Y$ is {\em horizontal} if and only if
\be
\ell(\gamma)=\ell(\pi\gamma).
\ee
The set of all such curves will be denoted by $\Hol(\pi)$.
\ed

\bd Given a curve $\alpha:[0,1]\rightarrow X$ and a point $u\in\pi\inv\alpha(0)$ a {\em parallel translation of $u$ along $\alpha$} is a horizontal $\gamma$ such that $\gamma(0)=u$ and $\pi\gamma=\alpha$.
\ed

\br It is easy to produce examples where given $\alpha$ and $u$ there exist no parallel translation as well as examples where there are even infinitely many such lifts.
\er

\bd Given a curve $\alpha:[0,1]\rightarrow X$ and a point $u\in\pi\inv\alpha(0)$ {\em the parallel translation of $u$ along $\alpha$} is given as relation $\Par\subseteq \pi\inv(\alpha(0))\times\pi\inv(\alpha(1))$. This can be regarded as a set-valued function 
$$ \Par^{\alpha}:\pi\inv(\alpha(0))\dashrightarrow\pi\inv(\alpha(1)),
$$
given by
\be \Par^{\alpha}(u)=\{\gamma(1)| \gamma\in\Hol(\pi), \pi\gamma=\alpha\}\subseteq\pi\inv(\alpha(1)).
\ee
\ed
\subsection{Algebraic considerations} Recall that in these singular spaces the natural generalization of the notion of parallelism is not given by well-defined functions but by relations between fibers over endpoints of curves in the base space.  

Even in this weak sense, composition and inverses of relations are well defined operations that endow the notion of parallelism (and even one of holonomy) certain richness in structure. 

The following definitions and statements are elementary (cf. \citet{MR1071176}) and are recalled here within the context of the parallel translation relations defined earlier in this section.

\bd Given sets $A$ and $B$, a {\em relation} $f:A\dashrightarrow B$ is a subset of $f\subseteq A\times B$. 
\ed

\br Looking at relations instead of at functions is sometimes advantageous; an example being the notion of correspondences given by Gromov in an equivalent formulation of the Gromov-Hausdorff distance. Usual functions $f:A\rightarrow B$  are of course particular cases of relations.
\er

\bd Given two relations $f:A\dashrightarrow B$, $g:B\dashrightarrow C$, the composition $g\circ f:A\dashrightarrow C$ is defined in the usual way:
\be
g\circ f=\{(a,c) | \exists b, (a,b)\in f, (b,c)\in g\}.
\ee
\ed

\bl The composition is associative and for any $f:A\dashrightarrow B$, 
\be
f\circ id_A=id_B\circ f= f
\ee
\el
\begin{proof} Elementary.
\end{proof}

\bd Given $f:A\dashrightarrow B$ there exists a relation ${f}^*:B\dashrightarrow A$ given by 
$$
{f}^*=\{ (b,a) | (a,b) \in f\}
$$
\ed

\br In the case of actual functions $f^*$ coincides with the inverse, whenever the latter exists.
\er

\bl Given two relations $f:A\dashrightarrow B$, $g:B\dashrightarrow C$, 
\be(g\circ f)^*={f}^*\circ{g}^*.
\ee
Furthermore, 
\be f^{**}=f.
\ee
\el
\begin{proof} Elementary.
\end{proof}
\bl Given relations $f,g,h,k$ such that $f\circ h$, $k\circ f$, $g\circ h$ and $k\circ h$ exist, then if 
\be
f\subseteq g
\ee 
then
\be
f\circ h\subseteq g\circ h
\ee
and
\be
k\circ f\subseteq k\circ h.
\ee

\el
\begin{proof}
Elementary.
\end{proof}
The following fact summarizes the previous statements. 
\bp\label{ordcatrel}
Given a set $X$, the set of relations, with the composition and subsumption given as before,  is an ordered ${}^*$-semigroup with identity.
\ep

\br It should be noted that this coincides with the usual holonomy group in the case of a metric of Sasaki-type, as well as in the case of Riemannian submersions in general (\cite{MR2500106}). 
\er

\br Given a relation $f:A\dashrightarrow B$, one can consider the set-valued function that assigns to each $a\in A$ the (possibly empty) of $b$ related to $a$. In the sequel, it will be denoted using the functional notation
\be
f(a)=\{b| (a,b)\in f\}.
\ee
\er

\bp A submonoid  $G$ of a relation monoid $M$ is a group if and only if for all $a\in G$, $a^*=a\inv$
\ep
\begin{proof}
The sufficiency is immediate since it prescribes the existence of an inverse, in particular it follows that for any $a\in G$, $a$ is an invertible function. For the necessity, one first sees that because $aa\inv=id$, then $a$ is necessarily onto, i.e. that for all $y$ there exists $x$, namely any element in $a\inv(y)$, such that $a(x)=y$. 

Furthermore, since the monoids are ordered (see Proposition \ref{ordcatrel}), if for $a,b\in G$ are such that $a\subseteq b$ then, by multiplication on both sides by $b\inv$ yields that
\be
ab\inv\subseteq id
\ee
which in turn implies equality since $ab\inv$ must be surjective. Therefore, since
\be
aa\inv,a\inv a\supseteq id,
\ee
equalities must hold as well.
\end{proof}

In the context of parallel translations one has the following fact.

\bt Given a submetry $\pi:Y\rightarrow X$, and given two curves $\alpha,\beta:I\rightarrow X$ such that $\alpha(1)=\beta(0)$, then
\be
\Par^{\beta\cdot\alpha}=\Par^{\beta}\circ\Par^{\alpha},
\ee
where $\beta\cdot\alpha$ stands for the concatenation of $\alpha$ and $\beta$. Also,
\be
\Par^{\alpha^-}=(\Par^{\alpha})^*,
\ee
where $\alpha^-$ is the reverse curve.  

Furthermore, given a fixed $x\in X$, the set
\be
\Hol_x:=\{\Par^{\alpha}|\alpha(0)=\alpha(1)=x\}
\ee
is a ${}^*$-semigroup with identity.
\et
\begin{proof} Because parallel translation is defined by horizontal curves, and the concatenation of curves is additive in length, it follows that the the concatenation of horizontal curves is horizontal, thus proving the first claim. The second claim follows by reversing the direction of the horizontal curves.

In particular, for the set of parallel translations along loops, since it is closed under composition and under the involution if follows that it is indeed a monoid.
\end{proof}

\bd\label{holmon} Given a submetry $\pi:Y\rightarrow X$ and a point $x\in X$, the monoid with involution 
\be
\Hol_x:=\{\Par^{\alpha}|\alpha(0)=\alpha(1)=x\}
\ee
will be called {\em Holonomy monoid} of $\pi$ at $x$.
\ed

\br The failure of these monoids from being groups will be looked at in the sequel. Also, this monoids will not be isomorphic in general. 
\er

\subsection{Norms, re-scalings and limits}
In this section, one further property inherited by the limits is analyzed, together with its implications to parallel translation. Namely, that of scalar multiplication. Fiberwise, this was already expected to happen, since individual fibers are obtained as a quotient of a Euclidean space by the linear action of the wane group. 

\bt For any sequence of Riemannian manifolds $\{(X_i,p_i)\}$ converging to $(X_{\infty},x_{\infty})$ consider a convergent family of bundles with metric connection $(E_i,h_i,\nabla_i)$ over it converging to $(E_{\infty},y_{\infty})$. There exists a continuous $\R$-action $$\R\times E_{\infty}\rightarrow E_{\infty}$$ such that there exists a subsequence of $\{E_i\}$ such that the standard $\R$-actions given by scalar multiplication converge uniformly on compact sets to it. 
\et
\begin{proof} The existence of said map follows from an application of the Arzel\`a-Ascoli theorem by the following reasoning.  Regarding $\R\times E_i$ as a metric space with the standard product metric, one sees that by requiring $a,b\in \R$, $u,v\in E_i$ such that 
$$
\sqrt{|a|^2+d^2(0_i,u)},\sqrt{|b|^2+d^2(0_i,v)}\leq R,
$$
for some fixed $R>>1$. Recall that the distance function on $E_i$ is given as in \eqref{distsaseq} and that therefore the distance between re-scalings of a common vector is bounded above by their linear distance, that is
\be
d(au,b,u)\leq\|u\|_i|a-b|.
\ee
Also, 
\begin{align}
d(bu,bv)&=\inf_{\alpha}\sqrt{\ell^2(\alpha)+|b|^2\|P^{\alpha}u-v\|^2}\nonumber\\
&=\max\{1,|b|\}\inf_{\alpha}\sqrt{\ell^2(\alpha)+\|P^{\alpha}u-v\|^2}\nonumber\\
&\leq Rd(u,v)\nonumber,
\end{align}

and therefore
\begin{align}
d(au,bv)&\leq d(au,bu)+d(bu,bv)\nonumber\\
&\leq R|a-b|+Rd(u,v)\nonumber\\
&\leq \sqrt{2}R \sqrt{|a-b|^2+d^2(u,v)}.\nonumber
\end{align}
This proves that the family of maps $(a,u)\mapsto au$ is equicontinuous when restricted to balls of a given radius. Thus by the Arzel\`a-Ascoli theorem, there exists a convergent subsequence and thus the required map exists. Furthermore, since for any $a$ the map $u\mapsto au$ is also a limit of the corresponding re-scaling maps, the defining properties of an $\R$-action are also verified, namely: For all $u\in E_{\infty}$ and for all $a,b\in\R$
\begin{align}
1\cdot u&=u, \\
a\cdot(b\cdot u)&=(ab)\cdot u 
\end{align}
\end{proof}

\br As expected, multiplication by zero is yields the {\em zero section} (as defined in Remark \ref{sasprops} and Proposition \ref{sascpt}), namely
\be
0\cdot u=\varsigma_{\infty}\circ\pi_{\infty}(u)
\ee
\er

\bc For any sequence of Riemannian manifolds $\{(X_i,p_i)\}$ converging to $(X_{\infty},x_{\infty})$ consider a convergent family of bundles with metric connection $(E_i,h_i,\nabla_i)$ over it converging to $(E_{\infty},y_{\infty})$. For any $u\in E_{\infty}$, the map $$t\mapsto tu,$$ for $t\geq0$ is a geodesic parametrized proportional to arc-length.
\ec
\begin{proof} This again follows from the fact that the corresponding maps for any sequence $\{u_i\}$, with $u_i\in E_i$ are geodesic rays (cf. Remark \ref{sasprops}) and an application of the Arzel\`a-Ascoli theorem. 
\end{proof}

\bc\label{absres} For any sequence of Riemannian manifolds $\{(X_i,p_i)\}$ converging to $(X_{\infty},x_{\infty})$ consider a convergent family of bundles with metric connection $(E_i,h_i,\nabla_i)$ over it converging to $(E_{\infty},y_{\infty})$, and let $\mu_{\infty}:E_{\infty}\rightarrow\R$ as in Proposition \ref{sascpt}. 

Then for any $u\in E_{\infty}$ and for any $a\in\R$,
\be
\mu_{\infty}(au)=|a|\mu_{\infty}(u)
\ee
\ec
\begin{proof} This can be verified in two ways: 1) Since $\mu_{\infty}$ is the limit of the norms and since scalar multiplication satisfies said equation at the level of norms, then so will the limit satisfy it; 2) In view of the previous corollary, since $\mu_{\infty}(v)$ is also the distance between $v$ and $0\cdot v=\varsigma_{\infty}\pi_{\infty}(v)$ and the map $t\mapsto tau$, being part of a geodesic ray, is a minimal geodesic.  
\end{proof}

\bc\label{cnstnorm} For any sequence of Riemannian manifolds $\{(X_i,p_i)\}$ converging to $(X_{\infty},x_{\infty})$ consider a convergent family of bundles with metric connection $(E_i,h_i,\nabla_i)$ over it converging to $(E_{\infty},y_{\infty})$. For any $u\in E_{\infty}$ there exists a sequence $\{u_i\}$, with $u_i\in E_i$ such that $\|u_i\|=\mu_{\infty}(u)$ for all $i$. 
\ec
\begin{proof}
There are two cases given by whether $\mu_{\infty}(u)=0$ or no. If it is then for any sequence of points $\{x_i\}$ converging to $\pi_{\infty}(u)$ it follows that $u_i=\varsigma_{\infty}(x_i)$ converges to $u$ as required.  If $\mu_{\infty}(u)\neq0$ then, without loss of generality, one can consider a sequence $\tilde{u_i}$ converging to $u$ such that for all $i$, $\|\tilde{u_i}\|\neq0$. Then by letting $u_i=(\mu_{\infty}(u)\left/\|\tilde{u_i}\|\right.)\tilde{u_i}$, the conclusion also follows since $\|\cdot\|$ converges to $\mu_{\infty}$.
\end{proof}

\bt Horizontal curves are the uniform limits of horizontal curves.
\et
\begin{proof} Let $\pi:E\rightarrow X$ be, as before, a limit of vector bundles $\pi_i:E_i\rightarrow X_i$ and let $\gamma$ be a horizontal curve between $u\in E$ and $v\in E$ and consider a sequence $\{\gamma_i\}$ of piecewise smooth curves converging uniformly to $\gamma$ such that their lengths $\ell(\gamma_i)$ converge to $\ell(\gamma)$ 
. Let $\alpha_i=\pi_i\gamma_i$, and let $\widetilde{\gamma_i}$ be the unique horizontal lifts of $\alpha_i$ with  $\widetilde{\gamma_i}(0)=\gamma_i(0)$. 

Because $\alpha_i=\pi_i\circ\gamma_i$ holds, it follows that $\alpha$ is the limit $\{\alpha_i\}$ and since $\widetilde{\gamma_i}$ have uniformly bounded lengths, one can assume that they converge uniformly.  Indeed, a uniform upper bound $C$ on lengths implies that 
$$\widetilde{\gamma_i}(t)\in B_C(\gamma_i(0))\subseteq E_i,$$
 and thus the convergence can be regarded as a Hausdorff convergence as in Remark \ref{GHisH}.

 Furthermore, it follows that if $\widetilde{\gamma}$ is the limit of $\widetilde{\gamma_i}$, then $\pi\widetilde{\gamma}=\alpha$. 

The claim is that $\widetilde{\gamma}=\gamma$. In fact, for each $i$, since the Riemannian structure on $\pi_i\inv\alpha_i$ is flat and Euclidean,  the curve 
\be
\vartheta_i:t\mapsto P^{\alpha_i}_t\left((1-t)\gamma_i(0)+t(P^{\alpha_i}_t)\inv(\gamma_i(1))\right)
\ee
is shorter than $\gamma_i$ with the same endpoints, and its length is given by
\be\label{endptccc}
\ell(\vartheta_i)=\sqrt{\ell^2(\alpha_i)+\|P^{\alpha}\left(\gamma_i(0)\right)-\gamma_i(1)\|^2}.
\ee

\begin{figure}[h!] 
\includegraphics[scale=0.7]{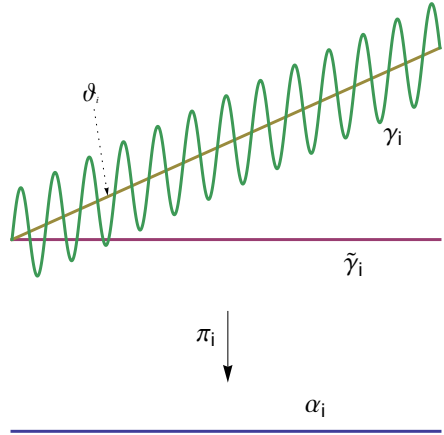}
\caption{The geometry on $\alpha_i^*E_i$}
\end{figure}

Now by the lower semi-continuity of the length functions, it follows that 
\be\ell(\gamma)=\lim_{i\rightarrow\infty}(\ell(\gamma_i))\geq\liminf_{i\rightarrow\infty}(\ell(\alpha_i))\geq \ell(\alpha)=\ell(\gamma),\ee
and that
\be
\ell(\gamma)=\lim_{i\rightarrow\infty}(\ell(\gamma_i))\geq\liminf_{i\rightarrow\infty}(\ell(\widetilde{\gamma_i}))\geq \ell(\widetilde{\gamma})\geq\ell(\alpha)=\ell(\gamma).
\ee
Therefore, $\widetilde{\gamma}$ is horizontal. For the curves $\vartheta_i$, since
\be
\ell(\gamma_i)\geq \ell(\vartheta_i)\geq\ell(\alpha_i),
\ee 
it follows from \eqref{endptccc} that
\be
\|P^{\alpha}(\gamma_i(0))-\gamma_i(1)\|\xrightarrow[i\rightarrow\infty]{} 0,
\ee
which means
\be
\|\tilde{\gamma}_i(1)-\gamma_i(1)\|\xrightarrow[i\rightarrow\infty]{} 0.
\ee
This now says that $\{\widetilde{\gamma_i}(1)\}$ converges to $\gamma(1)$, as required.  

Now, to see that this is true not only at the endpoints, notice that any segment of a horizontal curve is still horizontal, and that one can restrict the $\gamma_i$ and the $\alpha_i$ accordingly. The claim now follows and with it the end of the proof.
\end{proof}

\bc\label{absnorm} Horizontal curves have constant norm and constant re-scalings of horizontal curves are horizontal.
\ec
\begin{proof}
This is true since, before passing to the limit, the class of horizontal curves is closed under re-scaling and scalar multiplication is a uniform limit of scalar multiplications.
\end{proof}

\bt  For any sequence of Riemannian manifolds $\{(X_i,p_i)\}$ converging to $(X,p)$ consider a convergent family of bundles with metric connection $(E_i,h_i,\nabla_i)$ over it converging to $(E,\varsigma(p))$  with $\pi:E\rightarrow X$  as in Proposition \ref{sascpt}. Let $\alpha:I\rightarrow X$ be any rectifiable curve. Then for any $u\in \pi\inv(\alpha(0))$ there exists  $v\in \pi\inv(\alpha(1))$ such that there exists a horizontal path $\gamma$ from $u$ to $v$. In other words, 
\be
\Par^{\alpha}(u)\neq\varnothing.
\ee
\et
\begin{proof}
Let $\alpha_i$ be piecewise smooth curves converging uniformly to $\alpha$. Let $u_i\in\pi_i\inv(\alpha_i(0))$ be a sequence of points converging to $u\in \pi\inv(\alpha(0))$, Let $\gamma_i$ be the parallel translation along $\alpha_i$ with $\gamma_i(0)=u_i$. Since 
$$
\ell(\gamma_i)=\ell(\alpha_i)
$$
by construction, the convergence in length of $\{\alpha_i\}$ gives a uniform upper bound on the lengths of $\gamma_i$. Thus, by the Arzel\`a-Ascoli theorem, there exists a subsequence of $\{\gamma_i\}$, without loss of generality labeled again by $\gamma_i$, that converges uniformly to a curve $\gamma$ with $\gamma(0)=u$.

Now, by the lower semi-continuity of the the lengths,

\be
\ell(\alpha)=\lim{\ell(\alpha_i)}=\lim_{i\rightarrow\infty}\ell(\gamma_i)\geq\ell(\gamma)\geq\ell(\alpha),
\ee
which implies equality and thus finishes the proof.
\end{proof}

\subsection{Parallel translation along limit geodesics} In this section it will be seen that the limit of the total spaces of a sequence of vector bundles with corresponding Riemannian metrics of Sasaki-type contains many flats, i.e. isometric embeddings of $[a,b]\times[0,\infty)\subseteq\R^2$. 

\bd Let $X$ be a (pointed) limit of complete geodesic spaces $X_i$ (e.g. Riemannian manifolds)and $x,y\in X$. A curve $\alpha:[0,1]\rightarrow X$ is a {\em minimal limit geodesic} if 
\begin{enumerate}
\item $\alpha(0)=x$ and $\alpha(1)=y$; 
\item $\ell(\alpha)=d(x,y)$, i.e. if it is a minimal geodesic; and
\item there exist sequences $\{x_i\in X_i\}$, $\{y_i\in X_i\}$ and minimal geodesics $\alpha_i$, with $\alpha_i(0)=x_i$ and $\alpha_i(1)=y_i$ such that $\{\alpha_i\}$ converges to $\alpha$.
\end{enumerate}
\ed

\br Not every geodesic in $X$ is a limit geodesic. Limit geodesics have been studying by many (e.g.  \citet{MR1484888}).
\er


\bp \label{geodpar} For any sequence of Riemannian manifolds $\{(X_i,p_i)\}$ converging to $(X,p)$ consider a convergent family of bundles with metric connection $(E_i,h_i,\nabla_i)$ over it converging to $(E,\varsigma(p))$  with $\pi:E\rightarrow X$  as in Proposition \ref{sascpt}. Let $\alpha:I\rightarrow X$, parametrized by arc-length,  be a limit geodesic. Then for any $u\in \pi\inv(\alpha(0))$ there exists an isometric embedding 
\be
\varphi=\varphi_{u,\alpha}:[0,\ell(\alpha)]\times[0,\infty)\rightarrow E_{\infty}
\ee
with
\be
\pi_{\infty}\varphi(t,s)=\alpha(t)
\ee
and such that 
\be
u=\varphi(0,\mu_{\infty}(u)).
\ee
Furthermore, if $\gamma(t):=\varphi(t, \mu_{\infty}(u))$, then $\gamma$ is horizontal and
\be
\varphi(t,s)=\frac{s}{\mu_{\infty}(u)}\gamma(t)
\ee
\ep
\begin{proof} Let $\alpha_i$ be geodesics (parametrized by arc-length) such that $\alpha=\lim_{i\rightarrow\infty}\alpha_i$. Let $u_i\in\pi_i\inv(\alpha_i(0))$, with $\|u_i\|=\mu_i(u)$ and such that $\{u_i\}$ converges to $u$. Let $\gamma_i$ be the unique parallel translation along $\alpha_i$ with initial value $u_i$; that is that
\be
\gamma_i(t)=P_t^{\alpha_i}(u_i). 
\ee
It is only needed to show that the map
\be
\varphi_i(t,s)=\frac{s}{\|u_i\|}\gamma_i(t)
\ee
is an isometric embedding for each $i$. Once this is done an application of the Arzel\`a-Ascoli theorem completes the proof. 

Indeed, by \eqref{distsaseq} the distance between two images is given by 
\begin{align}
d(\varphi(t,s),\varphi(t',s'))&=\inf_{\beta}\sqrt{\ell^2(\beta)+\|P^{\beta}(\varphi(t,s))-\varphi(t',s')\|^2}\\
&=\inf_{\beta}\sqrt{\ell^2(\beta)+\frac{1}{\|u_i\|^2}\|(sP^{\beta}P_t^{\alpha}-s'P_{t'}^{\alpha})(u_i)\|^2},\label{mettt}\\
\end{align}
where the infimum is over curves $\beta$ connecting $\pi_i(\varphi(t,s))=\alpha_i(t)$ to $\pi_i(\varphi(t',s'))=\alpha_i(t')$.
Now, given that parallel translation is by linear isometries of the fiber, then 
\be
\|(sP^{\beta}P_t^{\alpha}-s'P_{t'}^{\alpha})(u_i)\|\geq |s-s'|\|u_i\|
\ee
since the closest points between the spheres of radius $s\|u_i\|$ and $s'\|u_i\|$ with common center is given by the right hand side. Also, for any $\beta$
\be
\ell(\beta)\geq |t-t'|
\ee
which is the distance between its endpoints. Thus
\be\label{melll}
d(\varphi(t,s),\varphi(t',s'))\geq\sqrt{|t-t'|^2+|s-s'|^2}.
\ee

However, by chosing $\beta=\alpha\left|[\min\{t,t'\},\max\{t,t'\}\right]$ in \eqref{mettt}, the reverse inequality from \eqref{melll} is obtained, thus yielding the claim and finishing the proof.
\end{proof}

\subsection{On the uniqueness of parallel translates} Since one can think of examples of non uniqueness (and in fact some are produced in \cite{GDSP}), it is only natural to wonder what conditions guarantee uniqueness in parallel translation.  In this section a necessary condition will be proved (essential  that all the fibers be homemorphic) and its equivalence with the holonomy monoids being groups.

In fact, the holonomy monoids defined in Definition \ref{holmon} already determine the non-uniqueness of parallel translation globally, as seen in the next result.

\bt\label{holmonpar}  Suppose $\pi:E\rightarrow X$ is a pointed Gromov-Hausdorff limit of a sequence metrics of Sasaki-type. The holonomy monoids are indeed groups if and only if parallel translation is unique.
\et

\begin{proof}Since, if at all, the inverse is given by ${}^*$, if follows that parallel translations along loops are functions if and only if holonomy monoids are groups. Now, if there is a curve $\alpha$ such that there are two distinct horizontal curves over it with same initial value but different endpoints, then $\alpha\inv\alpha$ will be a curve for which the parallel translation relation is necessarily a set-valued function. 
\end{proof}

The condition that parallel translations be unique already implies some further control on the possible collapses of the fibers, namely the following fact.

\bt\label{parwane} Suppose $\pi:E\rightarrow X$ is a pointed Gromov-Hausdorff limit of a sequence metrics of Sasaki-type of finitely many holonomy types.  For any $x\in X$ and consider its corresponding wane group $G_x\leq O(k)$, guaranteed by Theorem \ref{finholtype}. If parallel translations are unique then for all $x,y\in X$ their wane groups are equal,
\be G_x= G_y,
\ee 
up to conjugation by an element in $O(k)$.
\et
\begin{proof} 
Now, the contrapositive statement says that if the exist two points with non isomorphic groups then parallel translation is not unique. To prove this let $\{x_i\}, \{y_i\}\subseteq X_i$ be sequences that converge to $x,y\in X$ respectively, such that $G_x$ and $G_y$ are not isomorphic. 

Let $u_i\in\pi_i\inv(x_i)$ converging to some $u\in\pi\inv(x)$. To bring this to the level of holonomic spaces, fix a minimal geodesic $\alpha_i$ from $x_i$ to $y_i$ and isomorphisms

\begin{align}
\phi_i:\pi_i\inv(x_i)&\rightarrow \R^k\textup{,  and}\\
\varphi_i:\pi_i\inv(y_i)&\rightarrow \R^k,
\end{align}
such that for all $i,j$, 
\begin{align}
\phi_i(u_i)&=\phi_j(u_j) \textup{,  and} \\
 \varphi_i(P^{\alpha_i}u_i)&=\varphi_j(P^{\alpha_j}u_j). 
\end{align}
Let $\tilde u=\phi_i(u_i), v=\varphi_i(P^{\alpha_i}(u_i))$,  and 
\be
P_i=\varphi_i\circ P^{\alpha_i}\circ \phi_i\inv;
\ee 
thus for all $i$, 
\be
P_i(\tilde u)=v
\ee
Again without loss of generality,  $\{P_i\}$ converges in $O(k)$ to a map $P$, with 
\be
P(\tilde u)=v.
\ee

By assumption, since $G_x\ncong G_y$, there exists $g\in G_x$ such that  for all $h\in G_y$ (or reversely there exists $h\in G_y$ such that  for all $g\in G_x$),

\be
P(g \tilde u)\neq h P(\tilde u),
\ee

for otherwise $PG_yP\inv=G_x$, a contradiction to them being different. Suppose without loss of generality that it is not the parenthetical case, i.e. that it is $g\in G_x$ that exists. It follows that $g$ is necessarily not the identity map. 

Back at the level of fibers, this says that I can find elements ${\tilde u}_i=\phi_i\inv(\tilde u)\in \pi_i\inv(x_i)$ that converge to $u$ but such that their parallel translates along $\alpha_i$ remain away from the parallel translates along $\alpha_i$ of the $u_i$ by a definite amount. Passing to the limit (and taking a further subsequence if needed), the corresponding horizontal curves connecting them converge to distinct horizontal curves with the same starting point $u$.
\end{proof}
\bc Suppose $\pi:E\rightarrow X$ is a pointed Gromov-Hausdorff limit of a sequence metrics of Sasaki-type of finitely many holonomy types. If parallel translations are unique then all fibers of $\pi$ are homeomorphic.
\ec
\begin{proof}
This again follows from Theorem \ref{finholtype}, since for each $x\in X$ the topology of the fiber is determined by $G_x$.
\end{proof}

\bibliographystyle{plainnat}
\bibliography{ref}
\end{document}